\pdfoutput=1
\RequirePackage{silence}
\documentclass[a4paper,11pt,oneside]{amsart}
\usepackage[margin=2.5cm]{geometry}

\usepackage[utf8]{inputenc}

\usepackage{amssymb,amsmath,amsthm,amsfonts,mathrsfs,bbm}
\usepackage[table]{xcolor}
\usepackage{graphicx}
\usepackage{mathtools}
\usepackage{scalerel}
\usepackage{pifont}
\usepackage{booktabs}
\usepackage{tcolorbox}
\usepackage{stmaryrd}
\usepackage{arydshln}
\usepackage{colortbl}
\usepackage{caption}

\definecolor{spinach}{RGB}{46,139,87}
\definecolor{tomato}{RGB}{255,99,71}
\definecolor{pumpkin}{RGB}{224,180,80}

\definecolor{orchid}{RGB}{143,40,194}
\definecolor{lava}{RGB}{207,16,32}
\definecolor{mydarkblue}{RGB}{10,10,150}

\usepackage{xparse}

\usepackage{enumitem}
\setlist[enumerate]{itemsep=0.15cm,label=\emph{\upshape(\alph*).}}
\setlist[enumerate,2]{itemsep=0.15cm,label=\emph{\upshape(\roman*).}}

\setlist[description]{font=\normalfont, style=sameline,leftmargin=1cm,labelwidth=0.57cm}

\let\emph\relax
\DeclareTextFontCommand{\emph}{\em}

\usepackage{array}
\newcolumntype{C}{>{$}c<{$}}

\usepackage[vcentermath]{youngtab}

\usepackage[all]{xy}
\usepackage{tikz}

\definecolor{cream}{RGB}{255,253,208}

\usetikzlibrary{cd}
\usetikzlibrary{decorations}
\usetikzlibrary{decorations.markings}
\usetikzlibrary{decorations.pathreplacing}
\usetikzlibrary{decorations.pathmorphing}
\usetikzlibrary{arrows.meta,shapes,positioning,matrix,calc}
\usetikzlibrary{shapes.callouts}
\tikzset{anchorbase/.style={baseline={([yshift=-0.5ex]current bounding box.center)}},
tinynodes/.style={font=\tiny, text height=0.25ex, text depth=0.05ex},
smallnodes/.style={font=\scriptsize, text height=0.75ex, text depth=0.15ex},
mor/.style={line width=0.75,color=black,fill=cream},
crossline/.style={preaction={draw=white,line width=5.0pt,-},preaction={draw=black,line width=0.9pt,-}},
usual/.style={line width=1.0,color=black},
dot/.style = {
decoration={markings,
post length=0.25mm,
pre length=0.25mm,
mark=at position #1 with {\node[circle,radius=0.2cm,inner sep=-1.5pt,color=black,fill=black]{};}
},
postaction={decorate}
},
dot/.default=1,
}

\renewcommand{\dots}{\text{...}}

\newcommand{\C}{\mathbb{C}}
\newcommand{\Z}{\mathbb{Z}}
\newcommand{\R}{\mathbb{R}}
\newcommand{\Q}{\mathbb{Q}}
\newcommand{\F}{\mathbb{F}}
\newcommand{\K}{\mathbb{K}}
\newcommand{\LL}{\mathbb{L}}

\newcommand{\N}{\mathbb{Z}_{\geq 0}}
\newcommand{\cchar}[1][\K]{\mathrm{char}(#1)}

\newcommand{\End}{\mathrm{End}}
\newcommand{\Hom}{\mathrm{Hom}}

\newcommand{\monoid}[1][S]{\mathcal{#1}}

\newcommand{\onebt}{\mathbbm{1}_{bt}}

\newcommand{\gap}[2][\K]{\mathrm{gap}_{#1}(#2)}
\newcommand{\faith}[2][\K]{\mathrm{faith}_{#1}(#2)}

\newcommand{\gratio}[2][\K]{\mathrm{gapr}_{#1}(#2)}

\newcommand{\fratio}[2][\K]{\mathrm{faithr}_{#1}(#2)}

\newcommand{\jcell}{\mathcal{J}}
\newcommand{\hcell}{\mathcal{H}}

\newcommand{\rad}[1][\jcell]{\mathrm{Rad}}

\newcommand{\ssdimk}[1][\K]{\mathrm{ssdim}_{\K}}
\newcommand{\ssgap}[2][\K]{\mathrm{ssgap}_{#1}(#2)}
\newcommand{\ssratio}[2][\K]{\mathrm{ssgapr}_{#1}(#2)}

\newcommand{\tlalg}[2][n]{\mathcal{TL}_{#1}^{lin}(#2)}

\newcommand{\tlmon}[1][n]{\mathcal{TL}_{#1}}
\newcommand{\tltrun}{m\tlmon[n]^{k\le l}}

\newcommand{\brmon}[1][n]{\mathcal{B}\mathrm{r}_{#1}}

\newcommand{\pamon}[1][n]{\mathcal{P}\mathrm{a}_{#1}}
\newcommand{\robrmon}[1][n]{\mathcal{R}\mathrm{o}\mathcal{B}\mathrm{r}_{#1}}
\newcommand{\romon}[1][n]{\mathcal{R}\mathrm{o}_{#1}}
\newcommand{\ppamon}[1][n]{\mathrm{p}\mathcal{P}\mathrm{a}_{#1}}

\newcommand{\momon}[1][n]{\mathcal{M}\mathrm{o}_{#1}}
\newcommand{\promon}[1][n]{\mathrm{p}\mathcal{R}\mathrm{o}_{#1}}

\DeclarePairedDelimiterX{\inner}[2]{\langle}{\rangle}{#1, #2}

\DeclareMathOperator{\Id}{id}

\usepackage{aliascnt}
\def\NewTheorem#1{%
\newaliascnt{#1}{equation}%
\newtheorem{#1}[#1]{#1}%
\aliascntresetthe{#1}%
\expandafter\def\csname #1autorefname\endcsname{#1}%
}
\def\equationautorefname~#1\null{(#1)\null}

\numberwithin{equation}{subsection}

\NewTheorem{Proposition}
\NewTheorem{Theorem}
\newtheorem*{Theorem*}{Theorem}
\NewTheorem{Corollary}
\NewTheorem{Lemma}
\theoremstyle{definition}
\NewTheorem{Definition}
\NewTheorem{Notation}
\NewTheorem{Example}
\AtBeginEnvironment{Example}{%
  \pushQED{\qed}%
}
\AtEndEnvironment{Example}{\popQED}
\theoremstyle{remark}
\NewTheorem{Remark}
\AtBeginEnvironment{Remark}{%
  \pushQED{\qed}%
}
\AtEndEnvironment{Remark}{\popQED}
\NewTheorem{Assumption}
\NewTheorem{Hypothesis}
\NewTheorem{Conjecture}
\NewTheorem{Question}
\NewTheorem{Observation}
\NewTheorem{Fact}
\NewTheorem{Conclusion}

\usepackage[english]{babel}
\usepackage[hypertexnames=false]{hyperref}
\usepackage{bookmark}
\hypersetup{
pdftoolbar=true,
pdfmenubar=true,
pdffitwindow=false,
pdfstartview={FitH},
pdftitle={Representations of Cyclic Diagram Monoids},
pdfauthor={Jason Liu},
pdfsubject=,
pdfcreator={Jason Liu},
pdfproducer={Jason Liu},
pdfkeywords=,
pdfnewwindow=true,
colorlinks=true,
linkcolor=mydarkblue,
citecolor=teal,
filecolor=magenta,
urlcolor=orchid,
linkbordercolor=lava,
citebordercolor=teal,
urlbordercolor=orchid,
linktocpage=true
}

\addto\extrasenglish{%
}

\def\makeautorefname#1#2{\csdef{#1autorefname}{#2}}

\makeautorefname{section}{Section}%
\makeautorefname{subsection}{Section}%
\makeautorefname{subsubsection}{Section}%

\def\Today{\ifcase\month\or January\or February\or March\or
  April\or May\or June\or July\or August\or September\or
  October\or November\or December\fi\space\number\year}

\begin{document}



\title{Representations of Cyclic Diagram Monoids}
\author[J. Liu]{Jason Liu}

\address{Jason Liu: The University of Sydney, School of Mathematics and Statistics F07, NSW 2006, Australia}
\email{jliu0735@uni.sydney.edu.au}

\begin{abstract}
We introduce cyclic diagram monoids, a generalisation of classical diagram monoids that adds elements of arbitrary period by including internal components, with a view towards cryptography. We classify their simple representations and compute their dimensions in terms of the underlying diagram algebra. These go towards showing that cyclic diagram monoids possess representation gaps of exponential growth, which quantify their resistance as platforms against linear attacks on cryptographic protocols that exploit small dimensional representations.
\end{abstract}


\maketitle

\tableofcontents

\arrayrulewidth=0.5mm
\setlength{\arrayrulewidth}{0.5mm}


\section{Introduction}

Our main objective is to construct a generalisation of various diagram monoids that is potentially better suited for cryptography.

\subsection{Motivation}

One of the most prominent cryptographic protocols is the Diffie--Hellman (DH) key exchange, which enables two parties to establish a shared secret key. Its security is based on the difficulty of the discrete logarithm problem in finite fields, which underpins the security of many classical public-key protocols. However, the advent of quantum computing poses a significant threat to such schemes: Shor’s algorithm can solve the discrete logarithm problem in polynomial time, rendering protocols like the DH key exchange insecure in the presence of a sufficiently powerful quantum adversary. This has motivated the cryptographic community to search for alternative platforms that are resistant to quantum attacks.

Monoid and semigroup cryptography is studied as a way to explore new cryptographic primitives that might offer post-quantum security, non-commutative algebraic hardness, and novel computational problems that resist existing attacks. These structures allow for the design of protocols that do not rely on group inverses and often involve non-commutative operations, making them promising candidates for post-quantum cryptography. 

However, representation theory poses a potential backdoor attack on these protocols, which is known as a \emph{linear attack} (see e.g. \cite{romankov2014lineardecompositionattack}). Using a representation of the platform algebraic structure, the algorithmic problem underpinning the protocol can be solved using linear algebra and the solution is lifted back to the platform. For a linear attack to be efficient, the attacker needs a non-trivial representation of the monoid of small dimension (in a precise sense). Therefore, as a measurement of the security of the monoid as a platform, \cite{Khovanov_2024} defines the \emph{representation gap} of a monoid, which roughly measures the aforementioned dimension.

Diagram monoids provide interesting candidates with regards to representation gaps, which are families of monoids whose elements are diagrams with multiplication given by concatenation. The following table (from \cite{Khovanov_2024}) provides a summary of well-known diagram monoids in the literature.
\begin{gather*}\label{Eq:IntroDiaMonoids}
\begin{tabular}{c|c||c|c}
\arrayrulecolor{tomato}
Symbol & Diagrams 
& Symbol & Diagrams
\\
\hline
\hline
$\ppamon$ & \begin{tikzpicture}[anchorbase]
\draw[usual] (0.5,0) to[out=90,in=180] (1.25,0.45) to[out=0,in=90] (2,0);
\draw[usual] (0.5,0) to[out=90,in=180] (1,0.35) to[out=0,in=90] (1.5,0);
\draw[usual] (0.5,1) to[out=270,in=180] (1,0.55) to[out=0,in=270] (1.5,1);
\draw[usual] (1.5,1) to[out=270,in=180] (2,0.55) to[out=0,in=270] (2.5,1);
\draw[usual] (0,0) to (0,1);
\draw[usual] (2.5,0) to (2.5,1);
\draw[usual,dot] (1,0) to (1,0.2);
\draw[usual,dot] (1,1) to (1,0.8);
\draw[usual,dot] (2,1) to (2,0.8);
\end{tikzpicture}
& $\pamon$ & \begin{tikzpicture}[anchorbase]
\draw[usual] (0.5,0) to[out=90,in=180] (1.25,0.45) to[out=0,in=90] (2,0);
\draw[usual] (0.5,0) to[out=90,in=180] (1,0.35) to[out=0,in=90] (1.5,0);
\draw[usual] (0,1) to[out=270,in=180] (0.75,0.55) to[out=0,in=270] (1.5,1);
\draw[usual] (1.5,1) to[out=270,in=180] (2,0.55) to[out=0,in=270] (2.5,1);
\draw[usual] (0,0) to (0.5,1);
\draw[usual] (1,0) to (1,1);
\draw[usual] (2.5,0) to (2.5,1);
\draw[usual,dot] (2,1) to (2,0.8);
\end{tikzpicture}
\\
\hline
$\momon$ & \begin{tikzpicture}[anchorbase]
\draw[usual] (0.5,0) to[out=90,in=180] (1.25,0.5) to[out=0,in=90] (2,0);
\draw[usual] (1,0) to[out=90,in=180] (1.25,0.25) to[out=0,in=90] (1.5,0);
\draw[usual] (2,1) to[out=270,in=180] (2.25,0.75) to[out=0,in=270] (2.5,1);
\draw[usual] (0,0) to (1,1);
\draw[usual,dot] (2.5,0) to (2.5,0.2);
\draw[usual,dot] (0,1) to (0,0.8);
\draw[usual,dot] (0.5,1) to (0.5,0.8);
\draw[usual,dot] (1.5,1) to (1.5,0.8);
\end{tikzpicture}
& $\robrmon$ & \begin{tikzpicture}[anchorbase]
\draw[usual] (1,0) to[out=90,in=180] (1.25,0.25) to[out=0,in=90] (1.5,0);
\draw[usual] (1,1) to[out=270,in=180] (1.75,0.55) to[out=0,in=270] (2.5,1);
\draw[usual] (0,0) to (0.5,1);
\draw[usual] (2.5,0) to (2,1);
\draw[usual,dot] (0.5,0) to (0.5,0.2);
\draw[usual,dot] (2,0) to (2,0.2);
\draw[usual,dot] (0,1) to (0,0.8);
\draw[usual,dot] (1.5,1) to (1.5,0.8);
\end{tikzpicture}
\\
\hline
$\tlmon$ & \begin{tikzpicture}[anchorbase]
\draw[usual] (0.5,0) to[out=90,in=180] (1.25,0.5) to[out=0,in=90] (2,0);
\draw[usual] (1,0) to[out=90,in=180] (1.25,0.25) to[out=0,in=90] (1.5,0);
\draw[usual] (0,1) to[out=270,in=180] (0.25,0.75) to[out=0,in=270] (0.5,1);
\draw[usual] (2,1) to[out=270,in=180] (2.25,0.75) to[out=0,in=270] (2.5,1);
\draw[usual] (0,0) to (1,1);
\draw[usual] (2.5,0) to (1.5,1);
\end{tikzpicture}
& $\brmon$ & \begin{tikzpicture}[anchorbase]
\draw[usual] (0.5,0) to[out=90,in=180] (1.25,0.45) to[out=0,in=90] (2,0);
\draw[usual] (1,0) to[out=90,in=180] (1.25,0.25) to[out=0,in=90] (1.5,0);
\draw[usual] (0,1) to[out=270,in=180] (0.75,0.55) to[out=0,in=270] (1.5,1);
\draw[usual] (1,1) to[out=270,in=180] (1.75,0.55) to[out=0,in=270] (2.5,1);
\draw[usual] (0,0) to (0.5,1);
\draw[usual] (2.5,0) to (2,1);
\end{tikzpicture}
\\
\hline
$\promon$ & \begin{tikzpicture}[anchorbase]
\draw[usual] (0,0) to (0.5,1);
\draw[usual] (0.5,0) to (1,1);
\draw[usual] (2,0) to (1.5,1);
\draw[usual] (2.5,0) to (2.5,1);
\draw[usual,dot] (1,0) to (1,0.2);
\draw[usual,dot] (1.5,0) to (1.5,0.2);
\draw[usual,dot] (0,1) to (0,0.8);
\draw[usual,dot] (2,1) to (2,0.8);
\end{tikzpicture}
& $\romon$ & \begin{tikzpicture}[anchorbase]
\draw[usual] (0,0) to (1,1);
\draw[usual] (0.5,0) to (0,1);
\draw[usual] (2,0) to (2,1);
\draw[usual] (2.5,0) to (0.5,1);
\draw[usual,dot] (1,0) to (1,0.2);
\draw[usual,dot] (1.5,0) to (1.5,0.2);
\draw[usual,dot] (1.5,1) to (1.5,0.8);
\draw[usual,dot] (2.5,1) to (2.5,0.8);
\end{tikzpicture}
\end{tabular}
\end{gather*}

The diagram monoids on the left are planar versions of the ones on the right, i.e. their strands are not allowed to intersect. From top to bottom, we have: the partition monoids; the Motzkin and rook-Brauer monoids; the Temperley--Lieb and Brauer monoids; and the rook monoids. 

\cite{Khovanov_2024} demonstrated large representation gaps for many of these monoids. We are therefore motivated to study other aspects of their suitability as a platform monoid for cryptography. Some cryptographic protocols ask for elements of large periods (i.e. elements that generate a large cyclic group), but since not many well-known families of monoids contain elements of large periods, monoid-based protocols utilising periodicity are not well-studied. The elements of the planar diagram monoids above also all have period one, but these can be increased under a new construction, which can incentivise the development of new protocols utilising such elements. The periods of elements in a monoid are regulated by its \emph{Green's relations}, which partition a monoid into cells. By inflating these cells, we can construct variants of these diagram monoid that contain elements of high period. Therefore, we define \emph{cyclic diagram monoids}, in which internal components that appear in concatenation and are erased in the original diagram monoid are now counted modulo $m$, which inflates the aforementioned cells into cyclic groups. The goal is then to show that the cyclic version preserves all relevant properties of the original diagram monoid, including importantly the representation gap, which would justify their potential as platform monoids for cryptography.

\subsection{Results}

This paper is organised as follows. \autoref{Chapter1} recalls the necessary background about monoid representation theory,  representation gaps, and sandwich cellular algebras. Here we have a simple but useful new result in \autoref{T:FieldExtGap}. \autoref{Chapter3} is the focus of the paper, where we classify the simple representations of the cyclic Temperley--Lieb monoid and compute their dimensions in terms of the Temperley--Lieb algebra, thereby demonstrating the exponential growth of the representation gap of the (truncated) cyclic Temperley--Lieb monoid over fields of good characteristic, which generalises the result in \cite{Khovanov_2024} for the ordinary Temperley--Lieb monoid. \autoref{Chapter4} discusses the cyclic planar partition, Motzkin, and planar rook monoids, where we extend some of the same methods used for the cyclic Temperley--Lieb monoid, obtaining bounds for the representation gap of cyclic planar partition and rook monoids.

\subsection{Future directions}

There are more aspects to the suitability of cyclic diagram monoids for cryptography that need to be investigated. Importantly, their simple dimensions and representation gaps over fields of bad characteristic should be computed to demonstrate its security in all cases, which will require new approaches as, for example, the representation gap of the cyclic Temperley--Lieb monoid is not given by the dimension of the smallest non-trivial simple representation in this case. Additionally, we would like to be able to program the diagram monoids onto a computer with a computationally efficient multiplication, and   implement them into cryptographic protocols.

\subsection*{Acknowledgements}
Most of this paper was completed as part of the author's Master's thesis, who wishes to thank his supervisor, Dani Tubbenhauer, for their invaluable guidance and support throughout the project.

\section{Background}\label{Chapter1}

Let $\monoid$ be a monoid, and $\K$ be a field. For brevity, we have omitted proofs in this section, which can be found in \cite{Khovanov_2024} and \cite{Steinberg2016} for \autoref{S:Green} to \autoref{S:RepGapDef}, and \cite{Tubbenhauer_2024} for \autoref{S:Sandwich}.

\subsection{Green's relations}\label{S:Green}

Green's relations (see \cite{Green1951}) allow us to keep track of information loss during multiplication in a monoid.

\begin{Definition}
The \emph{left}, \emph{right}, and \emph{two-sided cell orders} on $\mathcal{S}$ are given by:
\begin{align*}
(a\leq_{l}b)&\Leftrightarrow
\exists c:b=ca,
\\
(a\leq_{r}b)&\Leftrightarrow
\exists c:b=ac,
\\
(a\leq_{lr}b)&\Leftrightarrow
\exists c,d:b=cad.
\end{align*}
\end{Definition}

We also write $<_l$, $\ge_r$, etc. to mean analogous orders.
(Note that we follow the convention used for cellular algebras, whereas in semigroup theory the orders are typically reversed.)

\begin{Definition}
 \emph{Green's relations} on $\mathcal{S}$  are the \emph{left}, \emph{right}, and \emph{two-sided equivalences} given by:
\begin{align*}
(a\sim_{l}b)&\Leftrightarrow
(a\leq_{l}b\text{ and }b\leq_{l}a),
\\
(a\sim_{r}b)&\Leftrightarrow
(a\leq_{r}b\text{ and }b\leq_{r}a),
\\
(a\sim_{lr}b)&\Leftrightarrow
(a\leq_{lr}b\text{ and }b\leq_{lr}a).
\end{align*}
An equivalence classes of the left (resp. right, two-sided) equivalence is called a left (resp. right, two-sided) \emph{cell}, which is denoted by $\mathcal{L}$ (resp. $\mathcal{R}$, $\mathcal{J}$). An $H$-cell is an intersection of a left cell and a right cell, denoted by $\mathcal{H}=\mathcal{H}(\mathcal{L},\mathcal{R})$.
\end{Definition}

In accordance with the notation, we also call two-sided cells $J$-cells. It also makes sense to write $a\le_l \mathcal{L}$ and $\mathcal{L}\le_l\mathcal{L'}$. We define $\mathcal{S}_{\le_\mathcal{L}}:=\{a\in \mathcal{S}\mid a\le_l\mathcal{L}\}$
and analogous sets by varying the subscript.

\begin{Definition}\label{D:Regular}
An element $e\in \mathcal{S}$ is an \emph{idempotent} if $e^2=e$. An $H$-cell or $J$-cell is \emph{idempotent} if it contains an idempotent. The monoid $\mathcal{S}$ is \emph{regular} if all $J$-cells of $\mathcal{S}$ are idempotent.
\end{Definition}

We denote idempotent $H$-cells and $J$-cells by $\mathcal{H}(e)$ and $\mathcal{J}(e)$ where $e$ is the idempotent it contains. 

\begin{Proposition}
The idempotent $H$-cell $\mathcal{H}(e)$ is a group with identity $ e$.
\end{Proposition}

\begin{Proposition}\label{L:UniqueJCell}
A monoid has a unique minimal $J$-cell and a unique maximal $J$-cell in the $\le_{lr}$-order, which are idempotent. In particular, the minimal $J$-cell is the group of units.
\end{Proposition}

We call the minimal $J$-cell the \emph{bottom cell} $\mathcal{J}_b$, and the maximal $J$-cell the \emph{top cell} $\mathcal{J}_t$. We also sometimes write $\mathcal{G}$ for $\mathcal{J}_b$ to emphasise that it is the group of units of $\mathcal{S}$.

Green's relations also provide information on the periodicity of elements, which is important in cryptography.

\begin{Definition}
Let $a\in\mathcal{S}$. The \emph{index} of $a$ is the smallest positive integer $i(a)$ such that $a^{i(a)}=a^{i(a)+d}$ for some $d\in\Z_{>0}$. The smallest possible such $d$ is the \emph{period} $p(a)$ of $a$.
\end{Definition}

Let $\mathcal{C}_m$ denote the cyclic group of order $m$.

\begin{Proposition}\label{T:Period}
Let $a\in\mathcal{S}$. There exists an idempotent $H$-cell $\mathcal{H}(e)$ such that $ \mathcal{C}_{p(a)}\cong \{ a^{s}\mid s\ge i(a) \}$ is a subgroup of $\mathcal{H}(e)$. In particular, $p(a)$ divides $|\mathcal{H}(e)|$.
\end{Proposition}

\subsection{Simple representations}\label{S:Simples}
We now recall the classification of simple representations for monoids. To begin, every monoid admits two trivial representations.

\begin{Definition}
The \emph{trivial representations} of $\mathcal{S}$ are:\begin{align*}\mathbbm{1}_{b}&:s\mapsto \begin{cases}1 & \text{if }s\in\mathcal{G}; \\0 & \text{else}.\end{cases} \\\mathbbm{1}_{t} & :s\mapsto 1\end{align*}
\end{Definition}

We write $\mathbbm{1}_{bt}$ to mean either $\mathbbm{1}_{b}$ or $\mathbbm{1}_{t}$, and $\mathbbm{1}_{bt}^{\oplus m}$ to mean any of the $2^m$ possible $m$-fold direct sums of $\mathbbm{1}_{bt}$.

\begin{Lemma}\label{L:SchutzenbergerReps}
For each left cell $\mathcal{L}$ of $\mathcal{S}$, there is a left $\mathcal{S}$-representation $\Delta_{\mathcal{L}}:=\K\mathcal{L}$ with $\K \mathcal{S}$-action:\[a\centerdot l=\begin{cases}al & \text{if }al\in \mathcal{L }; \\0 & \text{else}.\end{cases}\]
\end{Lemma}

These are known as \emph{Schützenberger representations}. 

\begin{Proposition}
Let $\mathcal{L}$ be a left cell of $\mathcal{S}$, and $\mathcal{H}(e)\subseteq\mathcal L$. Then the $\mathcal{S}$-representation $\Delta_{\mathcal{L}}$ is an $\mathcal{S}$-$\mathcal{H}(e)$-birepresentation.
\end{Proposition}

This allows us to construct $\mathcal{S}$-representations from $\mathcal{H}(e)$-representations.

\begin{Definition}
Let $M$ be an $\mathcal{H}(e)$-representation. The \emph{induced $\mathcal{S}$-representation} from $M$ is:
\[\operatorname{Ind}_{\mathcal{H}(e)}^{\mathcal{S}}M:=\Delta_{\mathcal{L}}\otimes_{\mathcal{H}(e)}M.\]
\end{Definition}

\begin{Definition}
Let $ M$ be an $\mathcal{S}$-representation. An \emph{apex} of $M$ is a maximal idempotent $ J$-cell $\mathcal{J}$ of $\mathcal{S}$ in the $\le_{lr}$-order that does not annihilate $M$.
\end{Definition}

\begin{Lemma}\label{L:UniqueApex}
Every simple $\mathcal{S}$-representation has a unique apex.
\end{Lemma}

Recall that the {head} $\operatorname{Hd}(M)$ is the quotient of $M$ by its radical. The following theorem classifies the simple $\mathcal{S}$-representations.

\begin{Theorem}[$H$-reduction]\label{T:CMP}
There is a bijection:\begin{align*}\{ \text{simple }\mathcal{H}(e)\text{-representations} \} &  \longleftrightarrow\{ \text{simple }\mathcal{S}\text{-representations of apex }\mathcal{J}(e) \} \\K & \,\longmapsto L_{K}:= \operatorname{Hd}(\operatorname{Ind}_{\mathcal{H}(e)}^{\mathcal{S}}K)\end{align*}
\end{Theorem}

An important structural property of monoids compared to groups is that they can be truncated in a way to keep only desired cells and hence only desired associated representations. 

\begin{Definition}
Let $\mathcal{J}$ be a $ J$-cell of $\mathcal{S}$ with $1\not\in \mathcal{J}$. The \emph{$\mathcal{J}$-submonoid} is $\mathcal{S}_{\ge\mathcal{J}}:=\mathcal{S}_{\ge_{lr}\mathcal{J}}\cup \{ 1 \}$.
\end{Definition}

Note that the $\mathcal{J}$-submonoid is indeed a submonoid of $\mathcal{S}$ since $\mathcal{S}_{\ge_{lr}\mathcal{J}}$ is a two-sided ideal.

\begin{Definition}
Let $ I$ be a two-sided ideal of $\mathcal{S}$. The \emph{Rees factor} $\mathcal{S}/I$ is the monoid $(\mathcal{S}\setminus I)\cup \{ 0 \}$ with multiplication:$$s\bullet t=\begin{cases}st & \text{if }st\in \mathcal{S}\setminus I; \\0 & \text{else.}\end{cases}$$
\end{Definition}

\begin{Definition}
Let $\mathcal{J}\le_{lr}\mathcal{K}$ be $ J$-cells of $\mathcal{S}$ with $1\not\in \mathcal{J}$. The \emph{$\mathcal{J}$-$\mathcal{K}$-subquotient} is the Rees factor $\mathcal{S_{J}^{K}}:=\mathcal{S_{\ge J}}/\mathcal{S_{\ge K}}$. We also define extremal cases $\mathcal{S}_{\mathcal{J}}:=\mathcal{S}_{\ge\mathcal{J}}$ and $ \mathcal{S}^{\mathcal{K}}:=\mathcal{S}/\mathcal{S}_{\ge\mathcal{K}}$. 
\end{Definition}

The following results can be stated similarly in the extremal cases.

\begin{Lemma}
Let $\mathcal{S}$ be regular, and $\mathcal{J}\le_{lr}\mathcal{K}$ be $ J$-cells of $\mathcal{S}$ with $1\not\in \mathcal{J}$. Then the $J$-cells of $\mathcal{S_{J}^{K}}$ are:
\[\{\mathcal{J}_b,\mathcal{J}_t\}\cup \{ \mathcal{M}\mid  \mathcal{M}\text{ a }J\text{-cell of }\mathcal{S}\text{ with }\mathcal{J}\le_{lr}\mathcal{M}<_{lr}\mathcal{K}\}.\]
Left, right, and $H$-cells are similarly classified.
\end{Lemma}

So when the monoid is regular, taking cell subquotients is essentially trimming away $J$-cells from the top and bottom of the monoid. This, along with the following proposition, motivates the use of monoids over groups for large representations, since we can keep only the part of the monoid that give large representations.

\begin{Proposition}\label{P:TruncatedRep}
Let $\mathcal{S}$ be regular, $\mathcal{J}\le_{lr}\mathcal{K}$ be $ J$-cells of $\mathcal{S}$ with $1\not\in \mathcal{J}$, and $\mathcal{M}\not\in \{ \mathcal{J}_{b},\mathcal{J}_{t} \}$ be an apex of $\mathcal{S_{J}^{K}}$. Then every simple $\mathcal{S_{K}^{J}}$-representation of apex $\mathcal{M}$ is a simple $\mathcal{S}$-representation of apex $\mathcal{M}$ (where $\mathcal{S}_{<_{lr}\mathcal{J}}$ acts by zero), and vice versa.
\end{Proposition}

\subsection{Representation gap}\label{S:RepGapDef}

Let $\mathcal{S}_{0;1}$ be the monoid $\{1,a\}$ with multiplication $a^2=a$, and recall that $\mathcal{C}_1$ denotes the trivial group.

\begin{Definition}
The pair $(\mathcal{S},\K)$, with $\mathcal{S}\not\cong\mathcal{C}_1,\mathcal{S}_{0;1}$, is \emph{$m$-trivial} if $\mathcal{S}$-representations $M$ with $\dim_\K(M)\le m$ satisfy $M\cong\onebt^{\oplus m}$. The \emph{representation  gap} $\operatorname{gap}_{\K}(\mathcal{S})$ of $(\mathcal{S},\K)$  is:
\[\operatorname{gap}_{\K}(\mathcal{S}):=\max\{m\mid (\mathcal{S},\K) \text{ is }(m-1)\text{-trivial}\}.\]
\end{Definition}

So a monoid with high representation gap will not admit any small dimensional non-trivial representation that can be exploited for a linear attack.

\begin{Remark}
The monoids $\mathcal{C}_1,\mathcal{S}_{0;1}$ are the only monoids whose representations are all of the form $\onebt^{\oplus n}$. We define them to be only $(-1)$-trivial, so that $\gap[\K]{\mathcal{C}_1}=\gap[\K]{\mathcal{S}_{0;1}}=0$.
\end{Remark}

We first summarise some simple but useful results that follow naturally from these definitions.

\begin{Lemma}\label{L:mTrivialEquiv}
The pair $(\mathcal{S},\K)$ is $ m$-trivial if and only if $ \mathcal{S}$-representations $ M$ with $\dim_{\K}(M)=m$ satisfy $M\cong\onebt^{\oplus m}$.
\end{Lemma}

\begin{Lemma}\label{L:GapBound}
Suppose $\mathcal{S}$ has a non-trivial simple representation. Then:
\[\operatorname{gap}_{\K}(\mathcal{S})\le \min \{ \dim_{\K}(L_{K})\mid L_{K}\not\cong \onebt\text{ a simple }\mathcal{S}\text{-representation} \}\le |\mathcal{S}|-1.\]
If $\K$ is algebraically closed, then $|\mathcal{S}|-1$ can be replaced by $\sqrt{|\mathcal{S}|-1}$. Moreover, if $\mathcal{S}$ is not a group, then every $|\mathcal{S}|-1$ can be replaced by ${|\mathcal{S}|-2}$.
\end{Lemma}

Another useful fact is that the representation gap over a field can be bounded below by the representation gap over a field extension. We first recall a classical theorem.

\begin{Theorem}[Noether--Deuring]\label{T:Noether}
Let $\mathscr{A}$ be a $\K$-algebra, $V,W$ be finite dimensional $\mathscr{A}$-representations, and $\LL$ be a field extension of $\K$. Then $V\otimes_\K \LL\cong W\otimes_\K \LL$ as $\mathscr{A}\otimes_\K\LL$-representations if and only if $V\cong W$ as $\mathscr{A}$-representations.
\end{Theorem}

\begin{Theorem}\label{T:FieldExtGap}
Let $\LL$ be a field extension of $\K$. Then: \[\gap[\LL]{\mathcal{S}}\le\gap[\K]{\mathcal{S}}.\]
\end{Theorem}
\begin{proof}
We show that $(\mathcal{S},\K)$ is $m$-trivial whenever $(\mathcal{S},\LL)$ is $m$-trivial, from which the result follows. Let $M$ be a $\K\mathcal{S}$-representation of dimension $m$. We have $\K\mathcal{S}\otimes_\K \LL\cong \LL\mathcal{S}$, and  $M\otimes_\K \LL$ is an $\LL\mathcal{S}$-representation also of dimension $m$. Suppose $(\mathcal{S},\LL)$ is $m$-trivial, then $M\otimes_\K \LL\cong \onebt^{\oplus m}\cong \onebt^{\oplus m}\otimes_\K \LL$ as $\LL\mathcal{S}$-representations. It then follows by the \hyperref[T:Noether]{Noether--Deuring theorem} that $M\cong \onebt^{\oplus m}$ as $\K\mathcal{S}$-representations, so $(\mathcal{S},\K)$ is $m$-trivial by \autoref{L:mTrivialEquiv}.
\end{proof}

Now recall the notion of cell subquotients.

\begin{Theorem}
Let $\mathcal{J}\le_{lr}\mathcal{K}$ be $J$-cells of $\mathcal{S}$. Then $\gap[\K]{\mathcal{S_J^K}}\ge \gap[\K]{\mathcal{S}}$.
\end{Theorem}

This result motivates the study of monoids over groups for large representation gaps: the presence of cell subquotients means we can truncate monoids to discard undesirable $J$-cells that are apexes to small representations. 

We now define more  measures of cryptographic security that are in the same vein as the representation gap.

\begin{Definition}
The \emph{faithfulness} of $(\mathcal{S},\K)$ is:
\[\operatorname{faith}_{\K}(\mathcal{S})=\min\{ \dim_{\K}(M)\mid M\text{ a faithful }\mathcal{S}\text{-representation} \}.\]
\end{Definition}

\begin{Remark}
Faithfulness has been studied as \emph{effective dimension} in \cite{mazorchuk2012effectivedimensionfinitesemigroups}. For groups, faithfulness over $\C$ has been studied more extensively as \emph{representation dimension} (see e.g. \cite{CerneleKamgarpourReichstein+2011+637+647}).
\end{Remark}

We have some useful bounds on faithfulness.

\begin{Lemma}\label{L:FaithBound}
Suppose $\mathcal{S}$ has a non-trivial simple representation. Then:
\[\gap[\K]{\mathcal{S}}\le\faith[\K]{\mathcal{S}}\le|\mathcal{S}|.\]
\end{Lemma}

\begin{Lemma}\label{L:FaithEmbed}
Let $\mathcal{S\hookrightarrow T}$ be an embedding of monoids. Then $\faith[\K]{\mathcal{S}}\le \faith[\K]{\mathcal{T}}$.
\end{Lemma}

In the case that the representation gap and faithfulness are difficult to find, we define another measure that is more crude but easier to compute. Recall the notation in \autoref{T:CMP}.

\begin{Definition}
Let $K$ be a simple $\mathcal{H}(e)$-representation. The \emph{semisimple dimension} of $L_K$ is: \[\ssdimk[\K](L_K)=\frac{|\mathcal{J}(e)|}{|\mathcal{H}(e)|}\dim_\K(K).\]
The \emph{semisimple representation gap} of $(\mathcal{S},\K)$ is:
\[\operatorname{ssgap}_{\K}(\mathcal{S})=\min\{ \ssdimk[\K](L_K)\mid L_K\not\cong\onebt\text{ a simple }\mathcal{S}\text{-representation} \}.\]
\end{Definition}

The following proposition justifies the name semisimple dimension.

\begin{Proposition}\label{P:ssdim}
The monoid algebra $\K\mathcal{S}$ is semisimple if and only if $\mathcal{S}$ is regular, all $\K\mathcal{H}(e)$ are semisimple, and $\dim_\K(L_K)=\ssdimk[\K](L_K)$ for all simple $\mathcal{S}$-representations $L_K$.
\end{Proposition}

We'd also like to take into account the size of the monoid as a measurement of efficiency.

\begin{Definition}
\begin{enumerate} 
    \item[]
    \item The \emph{gap ratio} of $\mathcal{S}$ is $\gratio[\K]{\mathcal{S}}:={\gap[\K]{\mathcal{S}}}/{\sqrt{|\mathcal{S}|}}$.
    \item The \emph{faithful ratio} of $\mathcal{S}$ is $\fratio[\K]{\mathcal{S}}:={\faith[\K]{\mathcal{S}}}/{|\mathcal{S}|}$.
    \item The \emph{semisimple gap ratio} of $\mathcal{S}$ is $\ssratio[\K]{\mathcal{S}}:={\ssgap[\K]{\mathcal{S}}}/{\sqrt{|\mathcal{S}|}}$.
\end{enumerate}
\end{Definition}

The square roots come from the case when $\K$ is algebraically closed. As a rule of thumb, a candidate family of monoids should not exceed polynomial decay in these ratios.

We have a sufficient condition for when the representation gap is well-behaved. Let $\Hom(\mathcal{S},\K)$ denote the set of monoid homomorphisms $\mathcal{S}\to(\K,+)$. In particular, there is always the zero map $0\in\Hom(\mathcal{S},\K)$. For $a,b\in\mathcal{S}\setminus \mathcal{G}$, we define a relation $ab\approx_{r}a$, and also write $\approx_{r}$ for its symmetric and transitive closure. We also define $\approx_{l}$ analogously.

\begin{Definition}
The monoid $\mathcal{S}$ is \emph{right-connected} if $\mathcal{S\setminus G}$ consists of a single equivalence class under $\approx_{r}$, \emph{left-connected} if $\mathcal{S\setminus G}$ consists of a single equivalence class under $\approx_{l}$, and \emph{null-connected} if every element of $\mathcal{S\setminus G}$ can be expressed as a product of two elements of $\mathcal{S\setminus G}$. We call $\mathcal{S}$ \emph{well-connected} if it is right, left, and null-connected, or if it is a group. 
\end{Definition}

\begin{Lemma}\label{L:RegNull}
If $\mathcal{S}$ is regular, then it is null-connected.
\end{Lemma}

\begin{Theorem}\label{T:WellConnect}
Let $\mathcal{S}$ be well-connected with $\Hom(\mathcal{S},\K)= \{0\}$. Then: \[\operatorname{gap}_{\K}(\mathcal{S})=\min\{ \dim_{\K}(L)\mid L\not\cong \onebt\text{ a simple }S\text{-representation} \}.\]
\end{Theorem}

\subsection{Sandwich cellular algebras}\label{S:Sandwich}
We assume some background on cellular algebras (see \cite{Graham1996}) and discuss a generalisation that is better suited for our purposes. Much of this follows a similar flavour as Green's relations on monoids, the reader might therefore like to skim the results for now and refer back to this subsection when they are used later on, as well as referring to \cite{Tubbenhauer_2024} for more detail.

\begin{Definition}

A $\K$-algebra $\mathscr{A}$ is \emph{sandwich cellular} if it is equipped with a sandwich cell datum $(\mathcal{P},(\mathcal{T,B}),(\mathscr{H}_{\lambda},B_{\lambda}),C)$ where:

\begin{enumerate}
\item $\mathcal{P}$ is the \emph{middle poset}.
\item $\mathcal{T}=\bigcup_{\lambda\in\mathcal{P}}\mathcal{T}(\lambda)$ and $\mathcal{B}=\bigcup_{\lambda\in\mathcal{P}}\mathcal{B}(\lambda)$ are collections of finite \emph{top} and \emph{bottom sets} respectively.
\item $\mathscr{H}_{\lambda}$ are the \emph{sandwiched algebras} with bases $B_{\lambda}$.
\item $C:\bigsqcup_{\lambda\in\mathcal{P}}\mathcal{T}(\lambda)\times B_{\lambda}\times \mathcal{B}(\lambda)\to \mathscr{A}$ given by $(T,m,B)\to c_{T,M,B}^{\lambda}$ is injective, and its image $B_{\mathscr{A}} $ is a $\K$-basis of $\mathscr{A}$ called the \emph{sandwich cellular basis}.
\end{enumerate}
such that:

\begin{enumerate}
\item For each $T\in\mathcal{T}(\lambda),\,m\in B_\lambda,\, B\in\mathcal{B}(\lambda)$, and $x\in\mathscr{A}$, we have:
\begin{equation}\label{SandwichCellularEq}
xc_{T,m,B}^{\lambda}\equiv \sum_{U\in\mathcal{T}(\lambda),\,n\in B_{\lambda}}r_{T,U}^{x}c_{U,n,B}^{\lambda}\pmod{\mathscr{A}^{>\lambda}},
\end{equation}
where $r_{T,U}^{x}\in\K$ do not depend on $B$ or $m$,  and $\mathscr{A}^{>\lambda}:=\K\{c_{T,m,B}^{\mu}\mid \mu>\lambda,\, T\in\mathcal{T}(\mu),\,m\in B_\mu,\, B\in\mathcal{B}(\mu)\}$. An analogous formula holds for right multiplication by $x$.
\item There exist a $\mathscr{A}$-$\mathscr{H}_{\lambda}$-birepresentation $\Delta(\lambda)$ and a $\mathscr{H}_{\lambda}$-$\mathscr{A}$-birepresentation $\nabla(\lambda)$ that are free with respect to $\mathscr{H}_\lambda$, such that there is an $\mathscr{A}$-$\mathscr{A}$-birepresentation isomorphism:\[\mathscr{A}_{\lambda}:={\mathscr{A}^{\ge\lambda}}/\mathscr{A}^{>\lambda}\cong \Delta(\lambda)\otimes_{\mathscr{H}_{\lambda}}\nabla (\lambda).\]
\end{enumerate}
\end{Definition}

\begin{Remark}
Multiplication in sandwich cellular algebras can be interpreted as diagrammatic concatenation by \autoref{SandwichCellularEq} in the following way:

\[\begin{tikzpicture}[anchorbase,scale=1]
\draw[mor] (0,-0.5) to (0.25,0) to (0.75,0) to (1,-0.5) to (0,-0.5);
\node at (0.5,-0.25){$B$};
\draw[mor] (0,1) to (0.25,0.5) to (0.75,0.5) to (1,1) to (0,1);
\node at (0.5,0.75){$T$};
\draw[mor] (0.25,0) to (0.25,0.5) to (0.75,0.5) to (0.75,0) to (0.25,0);
\node at (0.5,0.25){$m$};
\draw[mor] (0,1) to (0.25,1.5) to (0.75,1.5) to (1,1) to (0,1);
\node at (0.5,1.25){$B'$};
\draw[mor] (0,2.5) to (0.25,2) to (0.75,2) to (1,2.5) to (0,2.5);
\node at (0.5,2.25){$T'$};
\draw[mor] (0.25,1.5) to (0.25,2) to (0.75,2) to (0.75,1.5) to (0.25,1.5);
\node at (0.5,1.75){$m'$};
\end{tikzpicture}
\equiv
\underbrace{\begin{tikzpicture}[anchorbase,scale=1]
\draw[mor] (0,1) to (0.25,0.5) to (0.75,0.5) to (1,1) to (0,1);
\node at (0.5,0.75){$T$};
\draw[mor] (0,1) to (0.25,1.5) to (0.75,1.5) to (1,1) to (0,1);
\node at (0.5,1.25){$B'$};
\end{tikzpicture}}_{\in\K}
\cdot
\begin{tikzpicture}[anchorbase,scale=1]
\draw[mor] (0,-0.5) to (0.25,0) to (0.75,0) to (1,-0.5) to (0,-0.5);
\node at (0.5,-0.25){$B$};
\draw[mor] (0,1) to (0.25,0.5) to (0.75,0.5) to (1,1) to (0,1);
\node at (0.5,0.75){$T'$};
\draw[mor] (0.25,0) to (0.25,0.5) to (0.75,0.5) to (0.75,0) to (0.25,0);
\node at (0.5,0.25){\scalebox{0.55}{$m'm$}};
\end{tikzpicture}\pmod{\mathscr{A}^{>\lambda}}.\qedhere\]
\end{Remark}

We call $\mathscr{A}_\lambda$ a \emph{cell algebra}, and  $\Delta(\lambda),\nabla(\lambda)$ left and right \emph{cell representations}.

\begin{Definition}
A sandwich cellular algebra $\mathscr{A}$ is \emph{involutive} if it has an involutive sandwich cell datum $(\mathcal{P},\mathcal{T},(\mathscr{H}_{\lambda},B_{\lambda}),C,{*})$ where:

\begin{enumerate}
\item $\mathcal{T}(\lambda)=\mathcal{B}(\lambda)$ for all $\lambda\in \mathcal{P}$.
\item ${*}:\mathscr{A}\to \mathscr{A}$ is an anti-involution.
\item ${*}:B_{\lambda}\to B_{\lambda}$ is an bijection of order two.
\end{enumerate}
such that:$$(c_{T,m,B}^{\lambda})^{*}\equiv c_{B,m^{*},T}^{\lambda}\pmod{\mathscr{A}^{>\lambda}}.$$
\end{Definition}

From the definitions, we see that a cellular algebra is an involutive sandwich cellular algebra with $\mathscr{H}_\lambda\cong\K$ for all $\lambda\in\mathcal{P}$.

There is an analogue of Green's relations for algebras, such that we can reformulate sandwich cellularity in the language of cells. Let $\mathscr{A}$ be a finite dimensional algebra with a fixed basis $B_\mathscr{A}$. For $a,b,c\in B_\mathscr{A}$, we write $b\inplus ca$ if the coefficient of $b$ in the expansion of $ca$ with respect to $B_\mathscr{A}$ is non-zero.

\begin{Definition}\label{D:AlgebraCells}
The \emph{left}, \emph{right}, and \emph{two-sided cell orders} on $B_\mathscr{A}$ are given by:
\begin{align*}
(a\leq_{l}b)&\Leftrightarrow
\exists c:b\inplus ca,
\\
(a\leq_{r}b)&\Leftrightarrow
\exists c:b\inplus ac,
\\
(a\leq_{lr}b)&\Leftrightarrow
\exists c,d:b\inplus cad.
\end{align*}
The \emph{left}, \emph{right}, and \emph{two-sided equivalences} $\sim_l,\sim_r,\sim_{lr}$ on $B_\mathscr{A}$ are defined analogously to monoids, and so are \emph{left}, \emph{right}, \emph{two-sided}, and \emph{$H$-cells}.
\end{Definition}

\begin{Example}
Applying \autoref{D:AlgebraCells} to the monoid algebra $\K\mathcal{S}$ with respect to the basis $\mathcal{S}$ recovers Green's relations on $\mathcal{S}$.
\end{Example}

We define $\mathscr{A}^{>_{lr}\lambda}:=\K\{\bigcup_{\mu\in\mathcal{P},\,\mu>_{lr}\lambda}\mathcal{J}_\mu\}$, which is a two-sided ideal in $\mathscr{A}$. When the fixed basis $B_\mathscr{A}$ is a sandwich cellular basis of $\mathscr{A}$, we call the pair $(\mathscr{A},B_\mathscr{A})$ a \emph{sandwich pair}. If $B_\mathscr{A}$ admits an bijection of order two that induces an anti-involution on $\mathscr{A}$, then we call $(\mathscr{A},B_\mathscr{A})$ \emph{involutive}. By the analogy of Green's relations on both monoids and sandwich cellular algebras, we obtain a straightforward construction of sandwich cellular structures on monoid algebras.

\begin{Proposition}\label{P:MonSandPair}
Suppose $\mathcal{S}$ is regular or has singleton $H$-cells. Then $(\K\mathcal{S},\mathcal{S})$ is a sandwich pair. Furthermore, if $\mathcal{S}$ is anti-involutive, then $(\K\mathcal{S},\mathcal{S})$ is an involutive sandwich pair.
\end{Proposition}
\begin{proof}
We provide the construction for completeness. The middle poset $\mathcal{P}$ is the $\le_{lr}$-order from Green's relations on $\mathcal{S}$. The top sets $\mathcal{T}$ are indexed by right cells of $\mathcal{S}$, and the bottom sets $\mathcal{B}$ are indexed by left cells. The cell representations are $\Delta(\lambda)=\K\mathcal{L}$ and  $\nabla(\lambda)=\K\mathcal{R}$ for left and right cells $\mathcal{L,R}\subseteq\mathcal{J}_\lambda$. If $H$-cells are singletons, then $B_\lambda=\mathcal{H}\subseteq\mathcal{J}_\lambda$, else if $\mathcal{S}$ is regular, then $B_\lambda=\mathcal{H}(e)\subseteq\mathcal{J}_\lambda(e)$. Finally, the sandwich cellular basis consists of $c^\lambda_{T,m,B}=TmB$.
\end{proof}

Like cellular algebras, the multiplication in cell algebras of a sandwich cellular algebra is determined by a bilinear map on cell representations, which regulates its representation theory.

\begin{Lemma}
The multiplication in $\mathscr{A}_\lambda$ is determined by a bilinear map $\phi^\lambda: \nabla(\lambda)\otimes_{\mathscr{A}}\Delta(\lambda)\to\mathscr{H}_\lambda$.
\end{Lemma}

Diagrammatically, $\phi^\lambda$ is determined by concatenation as follows:
\[
\phi\left(
\scalebox{0.8}{$\begin{tikzpicture}[anchorbase,scale=1]
\draw[mor,orchid] (0,-0.5) to (0.25,0) to (0.75,0) to (1,-0.5) to (0,-0.5);
\node at (0.5,-0.25){$B^{\prime}$};
\draw[mor] (0,1) to (0.25,0.5) to (0.75,0.5) to (1,1) to (0,1);
\node at (0.5,0.75){$T$};
\draw[mor] (0.25,0) to (0.25,0.5) to (0.75,0.5) to (0.75,0) to (0.25,0);
\node at (0.5,0.25){$m$};
\end{tikzpicture}$}
,
\scalebox{0.8}{$\begin{tikzpicture}[anchorbase,scale=1]
\draw[mor] (0,-0.5) to (0.25,0) to (0.75,0) to (1,-0.5) to (0,-0.5);
\node at (0.5,-0.25){$B$};
\draw[mor,orchid] (0,1) to (0.25,0.5) to (0.75,0.5) to (1,1) to (0,1);
\node at (0.5,0.75){$T^{\prime}$};
\draw[mor] (0.25,0) to (0.25,0.5) to (0.75,0.5) to (0.75,0) to (0.25,0);
\node at (0.5,0.25){$m^{\prime}$};
\end{tikzpicture}$}
\right)
=
\scalebox{0.8}{$\begin{tikzpicture}[anchorbase,scale=1]
\draw[mor] (0,1) to (0.25,0.5) to (0.75,0.5) to (1,1) to (0,1);
\node at (0.5,0.75){$T$};
\draw[mor,orchid] (0,2.5) to (0.25,2) to (0.75,2) to (1,2.5) to (0,2.5);
\node at (0.5,2.25){$T^{\prime}$};
\draw[mor] (0.25,0) to (0.25,0.5) to (0.75,0.5) to (0.75,0) to (0.25,0);
\node at (0.5,0.25){$m$};
\draw[mor] (0,1) to (0.25,1.5) to (0.75,1.5) to (1,1) to (0,1);
\node at (0.5,1.25){$B$};
\draw[mor,orchid] (0,-0.5) to (0.25,0) to (0.75,0) to (1,-0.5) to (0,-0.5);
\node at (0.5,-0.25){$B^{\prime}$};
\draw[mor] (0.25,1.5) to (0.25,2) to (0.75,2) to (0.75,1.5) to (0.25,1.5);
\node at (0.5,1.75){$m^{\prime}$};
\end{tikzpicture}$}
\equiv
\underbrace{r_{TB}\cdot m^{\prime}m}_{\in\mathscr{H}_\lambda}
\pmod{\mathscr{A}^{>_{lr}\lambda}}
.
\]
The map $\phi^\lambda$ induces an $\mathscr{H}_\lambda$-linear map $\Delta(\lambda)\to\Hom_\K(\nabla(\lambda),\mathscr{H}_\lambda)$. Given an $\mathscr{H}_\lambda$-representation $K$, let $\Delta(\lambda,K):=\Delta(\lambda)\otimes_{\mathscr{H}_\lambda}K$. Then by tensoring with $\Id_K$, we can extend this map to: \[{\phi}^\lambda_K:\Delta(\lambda,K) \to\Hom_\K(\nabla(\lambda),\mathscr{H}_\lambda)\otimes_{\mathscr{H}_\lambda}K.\]
We denote the kernel of this map by $R(\lambda,K)$.

\begin{Definition}
Let $M$ be an $\mathscr{A}$-representation. An \emph{apex} of $M$ is a maximal $\lambda\in\mathcal{P}$ with $\K\mathcal{J}_\lambda$ not annihilating $M$.
\end{Definition}

\begin{Lemma}\label
Every simple $\mathscr{A}$-representation has a unique apex.
\end{Lemma}

\begin{Theorem}[$H$-reduction]\label{T:SandwichHReduce}
Let $\lambda\in\mathcal{P}$. If $\mathscr{H}_\lambda$ is finite dimensional, then there is a bijection:\begin{align*}\{ \text{simple }\mathscr{H}_\lambda\text{-representations} \} &  \longleftrightarrow\{ \text{simple }\mathscr{A}\text{-representations of apex }\lambda \} \\K & \,\longmapsto L(\lambda,K):=\Delta(\lambda,K)/R(\lambda,K)= \operatorname{Hd}(\Delta(\lambda,K))\end{align*}
\end{Theorem}

For the sandwich pair $(\K\mathcal{S},\mathcal{S})$, this recovers $H$-reduction for monoids.

\section{Cyclic Temperley--Lieb monoids}\label{Chapter3}

We will define the cyclic Temperley--Lieb monoids as endomorphism monoids of a family of monoidal categories, which follows the philosophy in \cite{Khovanov_2024} that interesting candidate families of monoids for cryptography arise from endomorphism monoids of monoidal categories.

\subsection{Definitions}\label{S:TLDefs}

Let $\equiv_m$ denote congruence modulo $m$ on $\Z$, and $\Z_m:=\Z/\mathord\equiv_m$. Recall that a compact connected 1-manifold is homeomorphic to either a closed interval or a circle. We now define the following diagrammatic category, which takes its name from \cite{HKAUFFMAN1987395}.

\begin{Definition}\label{Kauffman}
The \emph{Kauffman category} $\mathbf{K}$ is the monoidal category with:
\begin{enumerate}
\item Objects given by $\N$.

\item Morphisms $m\to n$ given by equivalence classes of \emph{crossingless matching} from $m$ points to $n$ points. A crossingless matching (or \emph{diagram} for short) is a compact 1-manifold embedded in $\R\times[0,1]$ with a finite number of connected components, such that its intersection with $\R\times\{1\}$ consists of $n$ \emph{top points} $\{(i,1)\mid i=1,\ldots,n)\}$ and its intersection with $\R\times\{0\}$ consists of $m$ \emph{bottom points} $\{(j,0)\mid j=1,\ldots,m)\}$. Two crossingless matchings $a,b$ are equivalent if there is a homeomorphism $h:a\to b$ fixing top and bottom points, i.e. $h(i,1)=(i,1)$ for all $i=1,\ldots,n$ and $h(j,0)=(j,0)$ for all $j=1,\ldots,m$.

\item Composition $\circ$ of $a:m\to n$ and $b: l\to m$ given by taking representative crossingless matchings and gluing along the $m$ middle points followed by compressing vertically.

\item Monoidal product $\otimes$ given on objects by $m\otimes n=m+n$, and on morphisms by taking representative crossingless matchings and placing one on the right of another (with top or bottom points shifted to the left accordingly).
\end{enumerate}
\end{Definition}

\begin{Example}
We can draw crossingless matchings as diagrams like so:
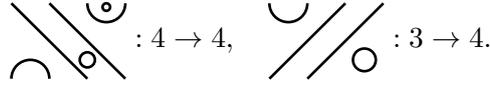
\[
\begin{tikzpicture}[anchorbase]
\draw[usual] (0,0) to[out=90,in=180] (0.25,0.25) to[out=0,in=90] (0.5,0);
\draw[usual] (1,1) to[out=270,in=180] (1.25,0.75) to[out=0,in=270] (1.5,1);
\draw[usual] (1,0) to (0,1);
\draw[usual] (1.5,0) to (0.5,1);
\draw[usual] (1,0.25) circle [radius=1mm];
\draw[usual] (1.25,0.95) circle [radius=0.5mm];
\end{tikzpicture}
:4\to 4,\quad 
\begin{tikzpicture}[anchorbase]
\draw[usual] (0,1) to[out=270,in=180] (0.25,0.75) to[out=0,in=270] (0.5,1);
\draw[usual] (0,0) to (1,1);
\draw[usual] (0.5,0) to (1.5,1);
\draw[usual] (1.25,0.25) circle [radius=1.5mm];
\end{tikzpicture}
:3\to 4.
\]

And their concatenation can be drawn as:
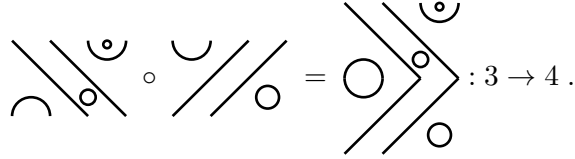
\[
\begin{tikzpicture}[anchorbase]
\draw[usual] (0,0) to[out=90,in=180] (0.25,0.25) to[out=0,in=90] (0.5,0);
\draw[usual] (1,1) to[out=270,in=180] (1.25,0.75) to[out=0,in=270] (1.5,1);
\draw[usual] (1,0) to (0,1);
\draw[usual] (1.5,0) to (0.5,1);
\draw[usual] (1,0.25) circle [radius=1mm];
\draw[usual] (1.25,0.95) circle [radius=0.5mm];
\end{tikzpicture}
\;\circ\; 
\begin{tikzpicture}[anchorbase]
\draw[usual] (0,1) to[out=270,in=180] (0.25,0.75) to[out=0,in=270] (0.5,1);
\draw[usual] (0,0) to (1,1);
\draw[usual] (0.5,0) to (1.5,1);
\draw[usual] (1.25,0.25) circle [radius=1.5mm];
\end{tikzpicture}
\;=\;
\begin{tikzpicture}[anchorbase]
\draw[usual] (0,1) to[out=90,in=180] (0.25,1.25) to[out=0,in=90] (0.5,1);
\draw[usual] (1,2) to[out=270,in=180] (1.25,1.75) to[out=0,in=270] (1.5,2);
\draw[usual] (1,1) to (0,2);
\draw[usual] (1.5,1) to (0.5,2);
\draw[usual] (1,1.25) circle [radius=1mm];
\draw[usual] (1.25,1.95) circle [radius=0.5mm];
\draw[usual] (0,1) to[out=270,in=180] (0.25,0.75) to[out=0,in=270] (0.5,1);
\draw[usual] (0,0) to (1,1);
\draw[usual] (0.5,0) to (1.5,1);
\draw[usual] (1.25,0.25) circle [radius=1.5mm];
\end{tikzpicture}:3\to 4\;.\qedhere
\]
\end{Example}

A crossingless matching has $\frac{n+m}{2}$ components homeomorphic to a closed interval with top and bottom points as endpoints, which we call \emph{strands}, and finitely many components homeomorphic to a circle, which we call \emph{internal components}. In particular, $\Hom_{\mathbf K}(m,n)=\emptyset$ if $m\not\equiv_2 n$.  We call a strand with both end points at the top a \emph{cup}, a strand with both end points at the bottom a \emph{cap}, and strand with one top and one bottom end point a \emph{through strand}.  Note that the identity morphism on $n$ is the crossingless matching with $n$ through strands and no internal components.

Another characterisation of equivalent crossingless matchings is if both diagrams have the same set of pairs of endpoints for their strands and the same number of internal components. Indeed, \autoref{Kauffman}.(b) implies:
\begin{gather*}
\begin{tikzpicture}[anchorbase]
\draw[usual] (0.5,0) to (0.5,1);
\draw[usual] (0.25,0.5) circle [radius=1mm];
\end{tikzpicture}
=
\begin{tikzpicture}[anchorbase]
\draw[usual] (0,0) to (0,1);
\draw[usual] (0.25,0.5) circle [radius=1mm];
\end{tikzpicture}
.
\end{gather*}
In particular, we can move internal components anywhere we want. If we let $C_m^n$ denote the set of equivalence classes of crossingless matchings from $m$ points to $n$ points but without internal components, then it follows that $\Hom_{\mathbf K}(m,n)$ is in bijection with $C_m^n\times \N$. We refer to the first component as the \emph{skeleton} of the crossingless matching. We define a map $\Phi$ on gluable sequences of crossingless matchings (i.e. consecutive diagrams having matching numbers of top and bottom points) that gives the number of internal components that would result from their concatenation. Then composition of morphisms in $\mathbf{K}$ transfers to the multiplication: \[(C_m^n\times \N)\times(C_l^m\times \N) \to (C_l^n\times \N),\quad  (a,i)(b,j)=(ab,i+j+\Phi(a,b)).\] 
Henceforth, we will be primarily using this formulation of the hom-sets. In addition, when we refer to a crossingless matching we will mean its equivalence class.

\begin{Lemma}\label{L:UniqueFactor}
Let $a\in \Hom_{\mathbf{K}}(m,n)$. Then there is a unique factorisation: \[a=(\gamma,0)(\Id_{k},l)(\beta,0)\] for $ k$ minimal, $(\beta,0)\in\Hom_{\mathbf{K}}(m,k)$, and $(\gamma,0)\in \Hom_{\mathbf K}(k,n)$.
\end{Lemma}
\begin{proof}
We take $k$ to be the number of through strands in $a$. The diagram $\gamma$ is given by the union of all cups and through strands, $\beta$ is given by the union of all caps and through strands, and $l$ is the number of internal components.  Their composition is then equivalent to $a$.
\end{proof}

Diagrammatically, the factorisation looks like:
\[
\begin{tikzpicture}[anchorbase]
\draw[usual] (1,0) to[out=90,in=180] (1.25,0.25) to[out=0,in=90] (1.5,0);
\draw[usual] (0.5,1) to[out=270,in=180] (0.75,0.75) to[out=0,in=270] (1,1);
\draw[usual] (-1,1) to[out=270,in=180] (-0.75,0.75) to[out=0,in=270] (-0.5,1);
\draw[usual] (0,0) to (0,1);
\draw[usual] (0.5,0) to (1.5,1);
\draw[usual] (-0.5,0.25) circle [radius=1.5mm];
\end{tikzpicture}
\;=\;
\underbrace{\begin{tikzpicture}[anchorbase]
\draw[usual] (0.5,1) to[out=270,in=180] (0.75,0.75) to[out=0,in=270] (1,1);
\draw[usual] (-1,1) to[out=270,in=180] (-0.75,0.75) to[out=0,in=270] (-0.5,1);
\draw[usual] (0,0) to (0,1);
\draw[usual] (1.5,0) to (1.5,1);
\end{tikzpicture}}_{(\gamma,0)}
\;\circ\;
\underbrace{\begin{tikzpicture}[anchorbase]
\draw[usual] (0,0) to (0,1);
\draw[usual] (0.5,0) to (0.5,1);
\draw[usual] (0.25,0.5) circle [radius=1mm];
\end{tikzpicture}}_{(\Id_{2},1)}
\;\circ\;
\underbrace{\begin{tikzpicture}[anchorbase]
\draw[usual] (1,0) to[out=90,in=180] (1.25,0.25) to[out=0,in=90] (1.5,0);
\draw[usual] (0,0) to (0,1);
\draw[usual] (0.5,0) to (0.5,1);
\end{tikzpicture}}_{(\beta,0)}
\;.
\]

We call $\gamma$ the \emph{top half} of the diagram, and $\beta$ the \emph{bottom half}. Note that $k$ is the number of through strands.

\begin{Definition}
The \emph{Kauffman monoid on $n$ strands}  is the endomorphism monoid: \[\mathcal{K}_n:=\End_{\mathbf{K}}(n)\cong C_n^n\times \N.\] 
\end{Definition}

\begin{Proposition}[{\cite[Lemma 3 \& 6]{borisavljevic2001kauffmanmonoids}}]\label{P:1stPres}
The monoid $\mathcal{K}_n$ has presentation:
\begin{equation*}
    \mathcal{K}_n= \scaleleftright[1.75ex]{<} {o,u_1,\ldots,u_{n-1}\, \vrule width 1pt \, \begin{matrix} u_i^2=ou_i=u_io,\\ u_iu_{i\pm 1}u_i=u_i,\\ u_iu_j=u_ju_i, \text{ for } |i-j|>1 \end{matrix}} {>}.
\end{equation*}
\end{Proposition}

Here $o$ corresponds to the diagram $(\Id_n,1)$ in the ordered pair notation. We call the $u_i$'s \emph{hooks}, which correspond to the following diagrams:
\[
\begin{tikzpicture}[anchorbase, scale=1.5, every node/.style={scale=0.75}]
\draw[usual] (0.5,1)node[above]{$i$} to[out=270,in=180] (0.75,0.75) to[out=0,in=270] (1,1)node[above]{$i+1$};
\draw[usual] (0.5,0) to[out=90,in=180] (0.75,0.25) to[out=0,in=90] (1,0);
\draw[usual] (0,0) to (0,1)node[above]{$i-1$};
\draw[usual] (1.5,0) to (1.5,1)node[above]{$i+2$};
\node at (-0.5,0.5) {$\cdots$};
\node at (2,0.5) {$\cdots$};
\draw[usual] (-1,0) to (-1,1)node[above]{$1$};
\draw[usual] (2.5,0) to (2.5,1)node[above]{$n$};
\end{tikzpicture}
\]

\begin{Definition}
The \emph{diagrammatic anti-involution} $*$ on $\mathbf{K}$ is the functor that sends a crossingless matching to its reflection across $\R\times\{\frac{1}{2}\}$.
\end{Definition}

This descends to an anti-involution $*$  on $\mathcal{K}_n$. The picture is:
\[
\left(\;\begin{tikzpicture}[anchorbase]
\draw[usual] (0,1) to[out=270,in=180] (0.25,0.75) to[out=0,in=270] (0.5,1);
\draw[usual] (0,0) to (1,1);
\draw[usual] (0.5,0) to (1.5,1);
\draw[usual] (1.25,0.25) circle [radius=1.5mm];
\end{tikzpicture}\;\right)^*
\;=\;
\begin{tikzpicture}[anchorbase, yscale=-1]
\draw[usual] (0,1) to[out=270,in=180] (0.25,0.75) to[out=0,in=270] (0.5,1);
\draw[usual] (0,0) to (1,1);
\draw[usual] (0.5,0) to (1.5,1);
\draw[usual] (1.25,0.25) circle [radius=1.5mm];
\end{tikzpicture}
\]

The Kauffman monoid is an infinite monoid, so we want to be able to take quotients to obtain finite monoids for cryptography. Recall that a congruence relation on a category is an equivalence relation on hom-sets that respects composition, and one can define the quotient category with the same objects and equivalence classes of morphisms as morphisms.

\begin{Lemma}
The binary relation $\simeq_m$ on $\Hom_{\mathbf K}(y,x)$ given by $(a,i)\simeq_m(a',i')$ if and only if $a=a'$ and $i\equiv_m i'$ is a congruence relation on $\mathbf{K}$. 
\end{Lemma}
\begin{proof}
We have an equivalence relation in each component, so it remains to check that $\simeq_m$ respects composition. Suppose $(a,i)\simeq_m(a',i')\in \Hom_{\mathbf K}(y,x)$ and $(b,j)\simeq_m(b',j')\in \Hom_{\mathbf K}(z,y)$, then:
\[(a,i)(b,j)=(ab,i+j+\Phi(a,b))\simeq_m (ab,i'+j'+\Phi(a,b))=(a,i')(b,j')=(a',i')(b',j'),\]
as needed.
\end{proof}

It follows that $\simeq_m$ is a congruence on the monoid $\mathcal{K}_n$ (in the sense of an equivalence relation that respects multiplication). In particular, we have the quotient monoid $\mathcal{K}_n/\mathord\simeq_m\cong C_n^n\times \Z_m$.

\begin{Definition}
The \emph{$m$-cyclic Temperley--Lieb category} is the quotient category: \[m\mathbf{TL}:=\mathbf{K}/\mathord\simeq_m.\]
\end{Definition}

By construction, $\simeq_m$ respects the monoidal product and diagrammatic involution on $\mathbf{K}$. Moreover, the unique factorisation of morphisms from \autoref{L:UniqueFactor} also descends to morphisms in $m\mathbf{TL}$. 

\begin{Lemma}\label{L:TLSize}
Let $\operatorname{Ca}(k)=\frac{1}{k+1}\binom{2k}{k}$ denote the $k$-th Catalan number. If $k\equiv_2 l$, then: \[|\Hom_{m\mathbf{TL}}(k,l)|=m\operatorname{Ca}(\frac{k+l}{2}).\]
\end{Lemma}
\begin{proof}
We have $o^m=1$ in the notation of \autoref{P:1stPres}, which implies that $o$ generates the cyclic group of order $m$. Now this follows from the well-known result that $|C_k^l|=\operatorname{Ca}(\frac{k+l}{2})$, see e.g. \cite[Section 5.7]{Kassel2008}.
\end{proof}

We can now define one of the central objects of our study. 

\begin{Definition}
The \emph{$m$-cyclic Temperley Lieb monoid on $n$ strands}  is the endomorphism monoid: \[m\tlmon[n]:=\End_{m\mathbf{TL}}(n)= \mathcal{K}_n/\mathord\simeq_m\cong C_n^n\times \Z_m.\]
\end{Definition}

In words, $m\tlmon[n]$ is the monoid of crossingless matchings with $n$ strands where the number of internal components is reduced modulo $m$. The case when $m=1$, i.e. internal components are removed outright, is denoted by $\tlmon[n]$.

\begin{Example}
Multiplication in $3\tlmon[4]$ looks like:
\[
\begin{tikzpicture}[anchorbase]
\draw[usual] (0,0) to[out=90,in=180] (0.25,0.25) to[out=0,in=90] (0.5,0);
\draw[usual] (1,1) to[out=270,in=180] (1.25,0.75) to[out=0,in=270] (1.5,1);
\draw[usual] (1,0) to (0,1);
\draw[usual] (1.5,0) to (0.5,1);
\draw[usual] (1,0.25) circle [radius=1mm];
\draw[usual] (1.25,0.95) circle [radius=0.5mm];
\end{tikzpicture}
\;\circ\; 
\begin{tikzpicture}[anchorbase]
\draw[usual] (0,1) to[out=270,in=180] (0.25,0.75) to[out=0,in=270] (0.5,1);
\draw[usual] (0,0) to (1,1);
\draw[usual] (0.5,0) to (1.5,1);
\draw[usual] (1.25,0.5) circle [radius=1mm];
\draw[usual] (1,0) to[out=90,in=180] (1.25,0.25) to[out=0,in=90] (1.5,0);
\end{tikzpicture}
\;=\;
\begin{tikzpicture}[anchorbase]
\draw[usual] (0,1) to[out=90,in=180] (0.25,1.25) to[out=0,in=90] (0.5,1);
\draw[usual] (1,2) to[out=270,in=180] (1.25,1.75) to[out=0,in=270] (1.5,2);
\draw[usual] (1,1) to (0,2);
\draw[usual] (1.5,1) to (0.5,2);
\draw[usual] (1,1.25) circle [radius=1mm];
\draw[usual] (1.25,1.95) circle [radius=0.5mm];
\draw[usual] (0,1) to[out=270,in=180] (0.25,0.75) to[out=0,in=270] (0.5,1);
\draw[usual] (0,0) to (1,1);
\draw[usual] (0.5,0) to (1.5,1);
\draw[usual] (1.25,0.5) circle [radius=1mm];
\draw[usual] (1,0) to[out=90,in=180] (1.25,0.25) to[out=0,in=90] (1.5,0);
\end{tikzpicture}
\;=\;
\begin{tikzpicture}[anchorbase]
\draw[usual] (1,0) to[out=90,in=180] (1.25,0.25) to[out=0,in=90] (1.5,0);
\draw[usual] (1,1) to[out=270,in=180] (1.25,0.75) to[out=0,in=270] (1.5,1);
\draw[usual] (0,0) to (0,1);
\draw[usual] (0.5,0) to (0.5,1);
\draw[usual] (1.25,0.5) circle [radius=1mm];
\end{tikzpicture}\;.
\]
where the last two diagrams are equivalent since they each have $4\equiv 1\pmod{3}$ internal components.
\end{Example}

\begin{Remark}\label{R:PartitionDiagram}
We can view diagrams in $m\tlmon[n]$ as partitions of $2n$ points $\{1,\ldots, n,1',\ldots,n'\}$ into sets of size 2, with strands connecting points in the same set and potential internal components.
\end{Remark}

\begin{Definition}
The \emph{Temperley--Lieb algebra} $\tlalg[n]{\delta}$  is the quotient algebra of $\K\mathcal{K}_n$ under the identification $o=\delta\cdot 1$ for some $\delta\in\K$.
\end{Definition}

In other words, in $\tlalg[n]{\delta}$  internal components are removed, and for each removed a multiplier of $\delta$ is introduced.

\begin{Remark}\label{R:TLMotivation}
The literature typically studies $\K$-linear versions of the Temperley--Lieb category with Temperley--Lieb algebras as hom-sets. However, algebras are inefficient as a platform for cryptography in terms of space complexity  due to their linearity. As a toy example, the monoid algebra over $\F_{p^n}$ has $p^n$ times the size of the monoid, but has representations over $\F_{p^n}$ of only the same dimensions. Indeed, \cite{Khovanov_2024} proposes producing set-theoretic counterparts of interesting algebras with large representations to fit the purpose of cryptography.
\end{Remark}

\begin{Remark}
In \cite{Coulembier_2023}, the growth rate of the number of indecomposable summands in tensor powers of a representation is studied, which motivates our study of cyclic Temperley--Lieb monoids as candidates for monoids with large representation gaps via the Schur--Weyl duality of $\tlalg[n]{\delta}\cong \End_{U_q(SL_2)}((\K^2)^{\otimes n})$, where $\delta=q+q^{-1}$ and $U_q(SL_2)$ is the quantum group of the special linear group $SL_2$. By e.g. \cite[Theorem 4.12]{Andersen_2018}, the dimensions of simple $\tlalg[n]{\delta}$-representations are given by the decomposition multiplicities of certain indecomposable $U_q(SL_2)$-representations in $(\K^2)^{\otimes n}$. Since \cite[Theorem 1.4]{Coulembier_2023} shows that the number of indecomposable summands grows exponentially, a priori some simple $\tlalg[n]{\delta}$-representations should be large and grow very fast. Since cyclic Temperley--Lieb monoids are constructed as set-theoretic counterparts of Temperley--Lieb algebras, we expect the same of them. It will also turn out that the large representations can be cut out of the monoid as per \autoref{S:Simples}.
\end{Remark}

\subsection{As platform monoid}\label{S:TLPlatform}

For cryptographic implementation, we would like to have a  presentation for $m\tlmon[n]$ in terms of generators and relations, and an efficient solution of the word problem with respect to the presentation. The former enables efficient computation of multiplications in the monoid, while the latter enables efficient common key establishment and verification. One way to solve the word problem is by having a normal form for elements of the monoid, which also helps disguise them from attackers. To this end, we consider an alternative presentation for $\mathcal{K}_n$.

\begin{Proposition}[{\cite[Section 1]{borisavljevic2001kauffmanmonoids}}]\label{P:2ndPres}
The monoid $\mathcal K_n$ is generated by $\{o\}\cup\{u^{[i,j]}\mid 1\le j\le i< n\}$ under the relations:
\begin{gather*}
u^{[i,j]}u^{[k,l]}=u^{[k,l]}u^{[i,j]},\text{ for }j\ge k+2, \\
u^{[i,j]}u^{[k,l]}=u^{[i,l]},\text{ for }i\ge l, |k-j|=1, \\
u^{[i,j]}o=ou^{[i,j]},\\
u^{[i,j]}u^{[j,l]}=ou^{[i,l]}.
\end{gather*}
\end{Proposition}

We refer to the $u^{[i,j]}$'s as \emph{blocks} and this presentation as the \emph{block formulation} of $\mathcal{K}_n$. The elements $u^{[i,j]}$ correspond to $u_iu_{i-1}\cdots u_{j+1}u_j$ in the hook formulation.

\begin{Proposition}[{\cite[Lemma 6]{borisavljevic2001kauffmanmonoids}}]\label{T:NormalForm}
Each element of $\mathcal{K}_n$ has a unique normal form $o^lu^{[b_1,a_1]}\cdots u^{[b_k,a_k]}$ with $a_1<\cdots<a_k$ and $b_1<\cdots<b_k$. In particular, $l$ is precisely the number of internal components in the corresponding diagram.
\end{Proposition}

The presentation and normal form descend to $m\tlmon[n]$ in the natural way. 

\begin{Theorem}\label{T:TLPres}
The monoid $m\tlmon[n]$ has the presentation of $\mathcal{K}_n$ with the additional relation that $o^m=1$.
\end{Theorem}
\begin{proof}
By \autoref{T:NormalForm}, the congruence $\simeq_m$ on the diagrammatic formulation of $\mathcal{K}_n$ passes to the congruence on the block formulation given by $o^lu^{[b_1,a_1]}\cdots u^{[b_k,a_k]}\simeq_m o^{l'}u^{[b'_1,a'_1]}\cdots u^{[b'_k,a'_k]}$ if and only if $l\equiv_m l'$ and $a_i=a_i',\,b_i=b_i'$ for all $1\le i\le k$. Consider the congruence closure of the relation $o^m=1$ on $\mathcal{K}_n$, i.e. the smallest congruence on $\mathcal{K}_n$ containing $o^m=1$. By induction, the congruence closure contains $o^{zm}=1$ for any $z\in\N$:
\[o^{zm}=1\Rightarrow o^{(z+1)m}=o^{zm}\cdot o^m=1\cdot 1=1.\]
So the congruence closure also contains:
\begin{align*}
o^lu^{[b_1,a_1]}\cdots u^{[b_k,a_k]}&=1\cdot(o^lu^{[b_1,a_1]}\cdots u^{[b_k,a_k]})\\&=o^{zm}\cdot(o^lu^{[b_1,a_1]}\cdots u^{[b_k,a_k]})\\&=o^{l+zm}u^{[b_1,a_1]}\cdots u^{[b_k,a_k]}.
\end{align*}
It follows that $o^m=1$ generates every relation in $\simeq_m$, so $\simeq_m$ is the congruence closure of $o^m=1$ on $\mathcal{K}_n$ as needed.
\end{proof}

So we obtain both a hook formulation and a block formulation  of $m\tlmon[n]$ from $\mathcal{K}_n$.

\begin{Corollary}
Each element of $m\tlmon[n]$ has a unique normal form $o^lu^{[b_1,a_1]}\cdots u^{[b_k,a_k]}$ with $0\le l<m$, $a_1<\cdots<a_k$, and $b_1<\cdots<b_k$.
\end{Corollary}
\begin{proof}
This follows from \autoref{T:NormalForm} and the proof of \autoref{T:TLPres} above.
\end{proof}

It follows that to solve the word problem in $m\tlmon[n]$, it suffices to reduce two words into their normal forms and compare them. Recall that the Big O notation $O(f)$ means asymptotically bounded above by $f$.

\begin{Proposition}
A word of length $l$ in $m\tlmon[n]$ can be reduced to its normal form in $O(l^2)$ operations.
\end{Proposition}
\begin{proof}
We can consider a word of $m\tlmon[n]$ as a word in $\mathcal{K}_n$ by viewing itself as a representative of its congruence class. By the proof of \cite[Lemma 1]{borisavljevic2001kauffmanmonoids}, every word in $\mathcal{K}_n$ can be reduced to its normal form by a set of rewriting rules that involves rewriting pairs of adjacent letters into either a single letter or another pair of letters. We can therefore look through the word from right to left until the first rewrite, then iterate until there are no more rewrites, and finally reduce the exponent of $o$ modulo $m$ to obtain the normal form in $m\tlmon[n]$. This is precisely the bubble sort algorithm (with the negligible speedup that the length of the word can decrease), which has $O(l^2)$ time complexity.
\end{proof}

\subsection{Cell structure}\label{S:TLCells}

We now work towards the representation gap of $m\tlmon[n]$, starting with some structural results. We begin with a simple lemma that reduces the proofs of many results to the $\tlmon[n]$ case. We identify elements of $\tlmon[n]$ with elements of $C_n^n$, that is, we refrain from writing them as ordered pairs.

\begin{Lemma}\label{L:BubOrder}
If $a\le_{l}b$ in $\mathcal{TL}_{n}$, then $(a,i)\le_{l}(b,j)$ in $m\mathcal{TL}_{n}$ for all $i,j\in\Z_m$. If $a\not\le_{l}b$ in $\mathcal{TL}_{n}$, then $(a,i)\not\le_{l}(b,j)$ in $m\mathcal{TL}_{n}$  for all $i,j\in\Z_m$. The same hold for $\le_{r}$ and $\le_{lr}$.
\end{Lemma}
\begin{proof}
Suppose $c  a=b$ in $\tlmon[n]$, then $(c,j-i-\Phi(c,a))(a,i)=(b,j)$ in $m\tlmon$ for all $i,j\in\Z_m$. Suppose $(c,k)(a,i)=(b,j)$ in $m\tlmon$ for some $i,j\in\Z_m$, then by definition $c  a= b$ in $\tlmon[n]$. The $\le_{r}$ and $\le_{lr}$ cases follow similarly.
\end{proof}

\begin{Corollary}\label{C:BubCell}
If $a\sim_{l}b$ in $\mathcal{TL}_{n}$, then $(a,i)\sim_{l}(b,j)$ in $m\mathcal{TL}_{n}$ for all $i,j\in\Z_m$. If $a\not\sim_{l}b$ in $\mathcal{TL}_{n}$, then $(a,i)\not\sim_{l}(b,j)$ in $m\mathcal{TL}_{n}$  for all $i,j\in\Z_m$. The same hold for $\sim_{r}$ and $\sim_{lr}$.
\end{Corollary}
\begin{proof}
This follows directly from \autoref{L:BubOrder}.
\end{proof}

It follows that $m\tlmon[n]$ has the cell structure of $\tlmon[n]$, but with each $H$-cell is inflated to size $m$. The cell structure is described in the following proposition. For $k\equiv_2 n$, we write $c(k):=(n-k)/2$.

\begin{Proposition}\label{P:BubCellStructure}
In $m\mathcal{TL}_{n}$:
\begin{enumerate}
    \item The left (resp. right) cells consist of crossingless matchings with the same bottom (resp. top) half. 

    \item The $J$-cells $\mathcal{J}_{k}$ consist of crossingless matchings with $k$ through strands.

        \item The $\le_{l}$, $\le_{r}$, and $\le_{lr}$-orders increase as the number of through strands decreases. In particular, $\le_{lr}$ is a total order.

    \item Within each $\mathcal{J}_{k}$, we have:
    \[|\mathcal{L}|=|\mathcal{R}|=m\left(\frac{k+1}{n-c(k)+1}\right)\binom{n}{c(k)}.\]

    \item For any $\mathcal{L}\subset \mathcal{J}_{k}$, we have $|\mathcal{J}_{k}|=|\mathcal{L}|^{2}/m$.

     \item Every $J$-cell is idempotent, i.e. $m\tlmon[n]$ is regular. 
     
     \item Every idempotent $H$-cell is isomorphic to $\mathcal{C}_m$.
\end{enumerate}
\end{Proposition}
\begin{proof}
\begin{description}
    \item[(a)--(e)] The $\tlmon[n]$ case is shown in e.g. \cite[Section 3]{Lau01082006}. The general result then follows from  \autoref{C:BubCell}.

    \item[(f)] We know $\tlmon[n]$ is regular by e.g. \cite[Section 3]{Lau01082006}, so let $a\in \tlmon[n]$ be an idempotent. Then: \[(a,-\Phi(a,a))^2=(a^2,-2\Phi(a,a)+\Phi(a,a))=(a,-\Phi(a,a).\] So $(a,-\Phi(a,a)$ is an idempotent.
    
     \item[(g)]  Let $(b,i)$ be an idempotent. From (a), we have $\mathcal{H}((b,i))=\{(b,0),\ldots, (b,m-1)\}$, so it remains to show it is cyclic. For all $j\in\Z_m$: \[(b,j)(b,1-\Phi(b,b))=(b,j+1-\Phi(b,b)+\Phi(b,b))=(b,j+1).\]  So $(b,1-\Phi(b,b))$ generates $\mathcal{H}((b,i))$. \qedhere
\end{description}
\end{proof}

Together, parts (f) and (g) of \autoref{P:BubCellStructure} ensure that $m\tlmon[n]$ contains elements of period $m$ as desired. Also, the proof of part (f) shows that $m\tlmon$ has the same idempotent $H$-cells as $\tlmon$.

\begin{Example}
The cell structure of $2\tlmon[4]$ is as follows:
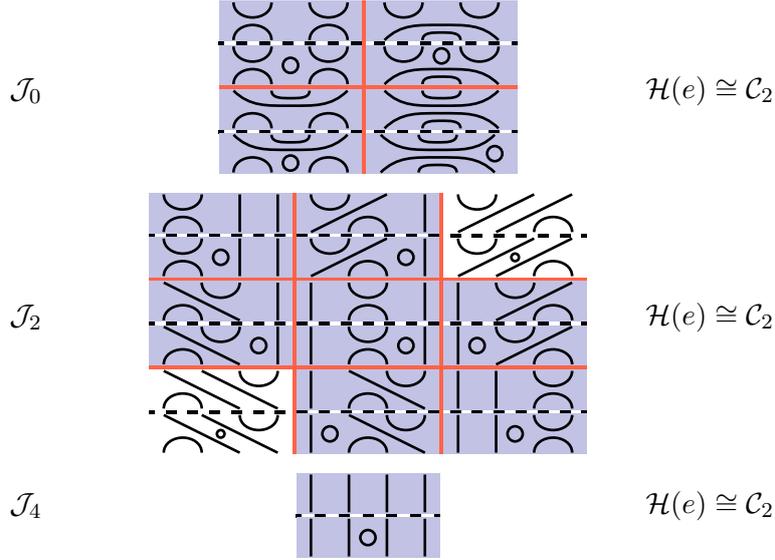
\begin{gather*}\label{Eq:TLCells}
\begin{gathered}
\xy
(0,0)*{\begin{gathered}\arrayrulecolor{tomato}
\begin{tabular}{C|C}
\cellcolor{mydarkblue!25}
\begin{tikzpicture}[anchorbase]
\draw[usual] (0,0) to[out=90,in=180] (0.25,0.2) to[out=0,in=90] (0.5,0);
\draw[usual] (0,0.5) to[out=270,in=180] (0.25,0.3) to[out=0,in=270] (0.5,0.5);
\draw[usual] (1,0) to[out=90,in=180] (1.25,0.2) to[out=0,in=90] (1.5,0);
\draw[usual] (1,0.5) to[out=270,in=180] (1.25,0.3) to[out=0,in=270] (1.5,0.5);
\end{tikzpicture} &
\cellcolor{mydarkblue!25}
\begin{tikzpicture}[anchorbase]
\draw[usual] (0,0) to[out=45,in=180] (0.75,0.20) to[out=0,in=135] (1.5,0);
\draw[usual] (0.5,0) to[out=90,in=180] (0.75,0.1) to[out=0,in=90] (1,0);
\draw[usual] (0,0.5) to[out=270,in=180] (0.25,0.3) to[out=0,in=270] (0.5,0.5);
\draw[usual] (1,0.5) to[out=270,in=180] (1.25,0.3) to[out=0,in=270] (1.5,0.5);
\end{tikzpicture}
\\
\arrayrulecolor{black}
\hdashline
\cellcolor{mydarkblue!25}
\begin{tikzpicture}[anchorbase]
\draw[usual] (0,0) to[out=90,in=180] (0.25,0.2) to[out=0,in=90] (0.5,0);
\draw[usual] (0,0.5) to[out=270,in=180] (0.25,0.3) to[out=0,in=270] (0.5,0.5);
\draw[usual] (1,0) to[out=90,in=180] (1.25,0.2) to[out=0,in=90] (1.5,0);
\draw[usual] (1,0.5) to[out=270,in=180] (1.25,0.3) to[out=0,in=270] (1.5,0.5);
\draw[usual] (0.75,0.25) circle [radius=1mm];
\end{tikzpicture} &
\cellcolor{mydarkblue!25}
\begin{tikzpicture}[anchorbase]
\draw[usual] (0,0) to[out=45,in=180] (0.75,0.20) to[out=0,in=135] (1.5,0);
\draw[usual] (0.5,0) to[out=90,in=180] (0.75,0.1) to[out=0,in=90] (1,0);
\draw[usual] (0,0.5) to[out=270,in=180] (0.25,0.3) to[out=0,in=270] (0.5,0.5);
\draw[usual] (1,0.5) to[out=270,in=180] (1.25,0.3) to[out=0,in=270] (1.5,0.5);
\draw[usual] (0.75,0.375) circle [radius=1mm];
\end{tikzpicture}
\\
\arrayrulecolor{tomato}
\hline
\cellcolor{mydarkblue!25}
\begin{tikzpicture}[anchorbase]
\draw[usual] (0,0) to[out=90,in=180] (0.25,0.2) to[out=0,in=90] (0.5,0);
\draw[usual] (1,0) to[out=90,in=180] (1.25,0.2) to[out=0,in=90] (1.5,0);
\draw[usual] (0,0.5) to[out=315,in=180] (0.75,0.3) to[out=0,in=225] (1.5,0.5);
\draw[usual] (0.5,0.5) to[out=270,in=180] (0.75,0.4) to[out=0,in=270] (1,0.5);
\end{tikzpicture} &
\cellcolor{mydarkblue!25}
\begin{tikzpicture}[anchorbase]
\draw[usual] (0,0) to[out=45,in=180] (0.75,0.20) to[out=0,in=135] (1.5,0);
\draw[usual] (0,0.5) to[out=315,in=180] (0.75,0.3) to[out=0,in=225] (1.5,0.5);
\draw[usual] (0.5,0) to[out=90,in=180] (0.75,0.1) to[out=0,in=90] (1,0);
\draw[usual] (0.5,0.5) to[out=270,in=180] (0.75,0.4) to[out=0,in=270] (1,0.5);
\end{tikzpicture}
\\
\arrayrulecolor{black}
\hdashline
\cellcolor{mydarkblue!25}
\begin{tikzpicture}[anchorbase]
\draw[usual] (0,0) to[out=90,in=180] (0.25,0.2) to[out=0,in=90] (0.5,0);
\draw[usual] (1,0) to[out=90,in=180] (1.25,0.2) to[out=0,in=90] (1.5,0);
\draw[usual] (0,0.5) to[out=315,in=180] (0.75,0.3) to[out=0,in=225] (1.5,0.5);
\draw[usual] (0.5,0.5) to[out=270,in=180] (0.75,0.4) to[out=0,in=270] (1,0.5);
\draw[usual] (0.75,0.125) circle [radius=1mm];
\end{tikzpicture} &
\cellcolor{mydarkblue!25}
\begin{tikzpicture}[anchorbase]
\draw[usual] (0,0) to[out=45,in=180] (0.75,0.20) to[out=0,in=135] (1.5,0);
\draw[usual] (0,0.5) to[out=315,in=180] (0.75,0.3) to[out=0,in=225] (1.5,0.5);
\draw[usual] (0.5,0) to[out=90,in=180] (0.75,0.1) to[out=0,in=90] (1,0);
\draw[usual] (0.5,0.5) to[out=270,in=180] (0.75,0.4) to[out=0,in=270] (1,0.5);
\draw[usual] (1.5,0.25) circle [radius=1mm];
\end{tikzpicture}
\end{tabular}
\\[3pt]
\begin{tabular}{C|C|C}
\arrayrulecolor{tomato}
\cellcolor{mydarkblue!25}
\begin{tikzpicture}[anchorbase]
\draw[usual] (0,0) to[out=90,in=180] (0.25,0.2) to[out=0,in=90] (0.5,0);
\draw[usual] (0,0.5) to[out=270,in=180] (0.25,0.3) to[out=0,in=270] (0.5,0.5);
\draw[usual] (1,0) to (1,0.5);
\draw[usual] (1.5,0) to (1.5,0.5);
\end{tikzpicture} & 
\cellcolor{mydarkblue!25}
\begin{tikzpicture}[anchorbase]
\draw[usual] (0.5,0) to[out=90,in=180] (0.75,0.2) to[out=0,in=90] (1,0);
\draw[usual] (0,0.5) to[out=270,in=180] (0.25,0.3) to[out=0,in=270] (0.5,0.5);
\draw[usual] (0,0) to (1,0.5);
\draw[usual] (1.5,0) to (1.5,0.5);
\end{tikzpicture} &
\begin{tikzpicture}[anchorbase]
\draw[usual] (1,0) to[out=90,in=180] (1.25,0.2) to[out=0,in=90] (1.5,0);
\draw[usual] (0,0.5) to[out=270,in=180] (0.25,0.3) to[out=0,in=270] (0.5,0.5);
\draw[usual] (0,0) to (1,0.5);
\draw[usual] (0.5,0) to (1.5,0.5);
\end{tikzpicture}
\\
\arrayrulecolor{black}
\hdashline
\cellcolor{mydarkblue!25}
\begin{tikzpicture}[anchorbase]
\draw[usual] (0,0) to[out=90,in=180] (0.25,0.2) to[out=0,in=90] (0.5,0);
\draw[usual] (0,0.5) to[out=270,in=180] (0.25,0.3) to[out=0,in=270] (0.5,0.5);
\draw[usual] (1,0) to (1,0.5);
\draw[usual] (1.5,0) to (1.5,0.5);
\draw[usual] (0.75,0.25) circle [radius=1mm];
\end{tikzpicture} & 
\cellcolor{mydarkblue!25}
\begin{tikzpicture}[anchorbase]
\draw[usual] (0.5,0) to[out=90,in=180] (0.75,0.2) to[out=0,in=90] (1,0);
\draw[usual] (0,0.5) to[out=270,in=180] (0.25,0.3) to[out=0,in=270] (0.5,0.5);
\draw[usual] (0,0) to (1,0.5);
\draw[usual] (1.5,0) to (1.5,0.5);
\draw[usual] (1.25,0.25) circle [radius=1mm];
\end{tikzpicture} &
\begin{tikzpicture}[anchorbase]
\draw[usual] (1,0) to[out=90,in=180] (1.25,0.2) to[out=0,in=90] (1.5,0);
\draw[usual] (0,0.5) to[out=270,in=180] (0.25,0.3) to[out=0,in=270] (0.5,0.5);
\draw[usual] (0,0) to (1,0.5);
\draw[usual] (0.5,0) to (1.5,0.5);
\draw[usual] (0.75,0.25) circle [radius=0.5mm];
\end{tikzpicture}
\\
\arrayrulecolor{tomato}
\hline
\cellcolor{mydarkblue!25}
\begin{tikzpicture}[anchorbase]
\draw[usual] (0,0) to[out=90,in=180] (0.25,0.2) to[out=0,in=90] (0.5,0);
\draw[usual] (0.5,0.5) to[out=270,in=180] (0.75,0.3) to[out=0,in=270] (1,0.5);
\draw[usual] (1,0) to (0,0.5);
\draw[usual] (1.5,0) to (1.5,0.5);
\end{tikzpicture} & 
\cellcolor{mydarkblue!25}
\begin{tikzpicture}[anchorbase]
\draw[usual] (0,0) to (0,0.5);
\draw[usual] (0.5,0) to[out=90,in=180] (0.75,0.2) to[out=0,in=90] (1,0);
\draw[usual] (0.5,0.5) to[out=270,in=180] (0.75,0.3) to[out=0,in=270] (1,0.5);
\draw[usual] (1.5,0) to (1.5,0.5);
\end{tikzpicture} &
\cellcolor{mydarkblue!25}
\begin{tikzpicture}[anchorbase]
\draw[usual] (0,0) to (0,0.5);
\draw[usual] (0.5,0) to (1.5,0.5);
\draw[usual] (1,0) to[out=90,in=180] (1.25,0.2) to[out=0,in=90] (1.5,0);
\draw[usual] (0.5,0.5) to[out=270,in=180] (0.75,0.3) to[out=0,in=270] (1,0.5);
\end{tikzpicture}
\\
\arrayrulecolor{black}
\hdashline
\cellcolor{mydarkblue!25}
\begin{tikzpicture}[anchorbase]
\draw[usual] (0,0) to[out=90,in=180] (0.25,0.2) to[out=0,in=90] (0.5,0);
\draw[usual] (0.5,0.5) to[out=270,in=180] (0.75,0.3) to[out=0,in=270] (1,0.5);
\draw[usual] (1,0) to (0,0.5);
\draw[usual] (1.5,0) to (1.5,0.5);
\draw[usual] (1.25,0.25) circle [radius=1mm];
\end{tikzpicture} & 
\cellcolor{mydarkblue!25}
\begin{tikzpicture}[anchorbase]
\draw[usual] (0,0) to (0,0.5);
\draw[usual] (0.5,0) to[out=90,in=180] (0.75,0.2) to[out=0,in=90] (1,0);
\draw[usual] (0.5,0.5) to[out=270,in=180] (0.75,0.3) to[out=0,in=270] (1,0.5);
\draw[usual] (1.5,0) to (1.5,0.5);
\draw[usual] (1.25,0.25) circle [radius=1mm];
\end{tikzpicture} &
\cellcolor{mydarkblue!25}
\begin{tikzpicture}[anchorbase]
\draw[usual] (0,0) to (0,0.5);
\draw[usual] (0.5,0) to (1.5,0.5);
\draw[usual] (1,0) to[out=90,in=180] (1.25,0.2) to[out=0,in=90] (1.5,0);
\draw[usual] (0.5,0.5) to[out=270,in=180] (0.75,0.3) to[out=0,in=270] (1,0.5);
\draw[usual] (0.25,0.25) circle [radius=1mm];
\end{tikzpicture}
\\
\arrayrulecolor{tomato}
\hline
\begin{tikzpicture}[anchorbase]
\draw[usual] (0,0) to[out=90,in=180] (0.25,0.2) to[out=0,in=90] (0.5,0);
\draw[usual] (1,0.5) to[out=270,in=180] (1.25,0.3) to[out=0,in=270] (1.5,0.5);
\draw[usual] (1,0) to (0,0.5);
\draw[usual] (1.5,0) to (0.5,0.5);
\end{tikzpicture} & 
\cellcolor{mydarkblue!25}
\begin{tikzpicture}[anchorbase]
\draw[usual] (0,0) to (0,0.5);
\draw[usual] (1.5,0) to (0.5,0.5);
\draw[usual] (0.5,0) to[out=90,in=180] (0.75,0.2) to[out=0,in=90] (1,0);
\draw[usual] (1,0.5) to[out=270,in=180] (1.25,0.3) to[out=0,in=270] (1.5,0.5);
\end{tikzpicture} &
\cellcolor{mydarkblue!25}
\begin{tikzpicture}[anchorbase]
\draw[usual] (0,0) to (0,0.5);
\draw[usual] (0.5,0) to (0.5,0.5);
\draw[usual] (1,0) to[out=90,in=180] (1.25,0.2) to[out=0,in=90] (1.5,0);
\draw[usual] (1,0.5) to[out=270,in=180] (1.25,0.3) to[out=0,in=270] (1.5,0.5);
\end{tikzpicture}
\\
\arrayrulecolor{black}
\hdashline
\begin{tikzpicture}[anchorbase]
\draw[usual] (0,0) to[out=90,in=180] (0.25,0.2) to[out=0,in=90] (0.5,0);
\draw[usual] (1,0.5) to[out=270,in=180] (1.25,0.3) to[out=0,in=270] (1.5,0.5);
\draw[usual] (1,0) to (0,0.5);
\draw[usual] (1.5,0) to (0.5,0.5);
\draw[usual] (0.75,0.25) circle [radius=0.5mm];
\end{tikzpicture} & 
\cellcolor{mydarkblue!25}
\begin{tikzpicture}[anchorbase]
\draw[usual] (0,0) to (0,0.5);
\draw[usual] (1.5,0) to (0.5,0.5);
\draw[usual] (0.5,0) to[out=90,in=180] (0.75,0.2) to[out=0,in=90] (1,0);
\draw[usual] (1,0.5) to[out=270,in=180] (1.25,0.3) to[out=0,in=270] (1.5,0.5);
\draw[usual] (0.25,0.25) circle [radius=1mm];
\end{tikzpicture} &
\cellcolor{mydarkblue!25}
\begin{tikzpicture}[anchorbase]
\draw[usual] (0,0) to (0,0.5);
\draw[usual] (0.5,0) to (0.5,0.5);
\draw[usual] (1,0) to[out=90,in=180] (1.25,0.2) to[out=0,in=90] (1.5,0);
\draw[usual] (1,0.5) to[out=270,in=180] (1.25,0.3) to[out=0,in=270] (1.5,0.5);
\draw[usual] (0.75,0.25) circle [radius=1mm];
\end{tikzpicture}
\end{tabular}
\\[3pt]
\begin{tabular}{C}
\cellcolor{mydarkblue!25}
\begin{tikzpicture}[anchorbase]
\draw[usual] (0,0) to (0,0.5);
\draw[usual] (0.5,0) to (0.5,0.5);
\draw[usual] (1,0) to (1,0.5);
\draw[usual] (1.5,0) to (1.5,0.5);
\end{tikzpicture}
\\
\arrayrulecolor{black}
\hdashline
\cellcolor{mydarkblue!25}
\begin{tikzpicture}[anchorbase]
\draw[usual] (0,0) to (0,0.5);
\draw[usual] (0.5,0) to (0.5,0.5);
\draw[usual] (1,0) to (1,0.5);
\draw[usual] (1.5,0) to (1.5,0.5);
\draw[usual] (0.75,0.25) circle [radius=1mm];
\end{tikzpicture}
\end{tabular}
\end{gathered}};
(-45,25)*{\jcell_{0}};
(-45,-5)*{\jcell_{2}};
(-45,-30)*{\jcell_{4}};
(45,25)*{\hcell(e)\cong\mathcal{C}_2};
(45,-5)*{\hcell(e)\cong\mathcal{C}_2};
(45,-30)*{\hcell(e)\cong\mathcal{C}_2};
(-51,0)*{\phantom{a}};
\endxy
\quad
\end{gathered}
\end{gather*}
where the dashed borders separate elements, red borders separate the $H$-cells, and idempotent $H$-cells are coloured. The columns are left cells and the rows are right cells.
\end{Example}

\begin{Proposition}\label{P:No.Reps}
The set of apexes for simple $ m\tlmon[n]$-representations can be indexed by the poset $\{ n,n-2,\dots \}$ under the reverse natural order. Moreover, there are at most $m$ simple $m\mathcal{TL}_{n}$-representation of a fixed apex. In particular, there are $m$ simple $m\mathcal{TL}_{n}$-representation of a fixed apex over a field that contains a primitive $m$-th root of unity.
\end{Proposition}
\begin{proof}
By \hyperref[T:CMP]{$H$-reduction} and \autoref{P:BubCellStructure}, the simple $m\mathcal{TL}_{n}$-representations of a fixed apex are in bijection with simple $\mathcal C_m$-representations. By standard representation theory, a group has at most as many simple representations as its order.
\end{proof}

Recall the notation for Schützenberger representations from \autoref{L:SchutzenbergerReps}.

\begin{Lemma}
Within the same $J$-cell, all $\Delta_{\mathcal{L}}$ (resp. $\Delta_{\mathcal{R}}$) are isomorphic.
\end{Lemma}
\begin{proof}
The action of $m\tlmon[n]$ on $\Delta_{\mathcal{L}}$ is by concatenation on top, which only depends on the top half of a diagram. \autoref{P:BubCellStructure} says that each left cell consists of diagrams with the same bottom half, so $m\tlmon[n]$ has the same action on each $\Delta_{\mathcal{L}}$. The same argument applies to right cells.
\end{proof}

We write $\Delta_{k}$ (resp. $_{k}\Delta$) for an arbitrary $\Delta_{\mathcal{L}}$ (resp. $\Delta_{\mathcal{R}}$) in $\mathcal{J}_{k}$.

\begin{Lemma}
We have:
\[\dim_{\K}(\Delta_{k})=\dim_{\K}(_{k}\Delta)=m\left(\frac{k+1}{n-c(k)+1}\right)\binom{n}{c(k)}.\]
\end{Lemma}
\begin{proof}
This follows from \autoref{L:SchutzenbergerReps} and \autoref{P:BubCellStructure}.
\end{proof}

\begin{Proposition}\label{P:SSDimFormula}
For $\mathcal{H}(e)\subseteq \mathcal{J}_{k}$ and $K$ a simple $\mathcal{H}(e)$-representation:
\[(m-1)\frac{k+1}{n-c(k)+1}\binom{n}{c(k)}\ge \operatorname{ssdim}_{\K}(L_{K})\ge \frac{k+1}{n-c(k)+1}\binom{n}{c(k)}.\]
In particular, if $\K$ contains a primitive $m$-th root of unity, then:
\[\operatorname{ssdim}_{\K}(L_{K})= \frac{k+1}{n-c(k)+1}\binom{n}{c(k)}.\]
\end{Proposition}
\begin{proof}
Both equations follow from \autoref{P:BubCellStructure} and \autoref{P:No.Reps}.
\end{proof}

\begin{Remark}
The first inequality above is sharp since $\mathcal{C}_p$ with $p$ prime has a simple representation of dimension $p-1$ over $\Q$.
\end{Remark}

If $\K$ contains an $m$-th root of unity, then $m\tlmon[n]$ has a non-trivial one-dimensional representation, namely the one where the generator of the group of units acts by multiplication by that root of unity and the non-invertible elements act by zero, which has apex $\mathcal{J}_b$. We'd like to truncate $m\tlmon[n]$ into cell subquotients using  \autoref{P:TruncatedRep} to discard these undesirable representations of small dimension (which we can do since $m\tlmon[n]$ is regular by \autoref{P:BubCellStructure}). Recall the notation for cell subquotients that $ \mathcal{S}^{\mathcal{K}}_{\mathcal J}=\mathcal{S}_{\ge\mathcal{J}}/\mathcal{S}_{\ge\mathcal{K}}$.

\begin{Definition}
For $0\le k< l\le  n$ with $k,l\equiv_2 n$. The \emph{$(k,l)$-th truncated Temperley--Lieb monoid} is the cell subquotient:
\begin{equation*}
    m\tlmon[n]^{k\le l}:=(m\tlmon[n])_{\mathcal{J}_l}^{\mathcal{J}_{k-2}}.
\end{equation*}
\end{Definition}

In words, $m\tlmon[n]^{k\le l}$ is the cell subquotient consisting of diagrams with between $k$ and $l$ through strands, with any product containing less than $k$ through strands killed. When $k=0,1$, we will write $m\tlmon[n]^{\le l}:=m\tlmon[n]^{k\le l}$ as the truncated monoid is the cell submonoid $(m\tlmon[n])_{\mathcal{J}_l}$ in this case.

\subsection{Trivial extensions}\label{S:TLExtension}

We now show that $\tltrun$ satisfies the connectivity conditions of \autoref{T:WellConnect}.  To begin, null-connectedness is immediate.

\begin{Lemma}\label{L:TLNullConn}
The monoid $\tltrun$ is null-connected.
\end{Lemma}
\begin{proof}
This follows from \autoref{L:RegNull} and \autoref{P:BubCellStructure}. 
\end{proof}

We now work towards left and right-connectedness.

\begin{Lemma}\label{L:DiffBubLeftConn}
In $\tltrun$, we have $(a,i)\approx_{l}(a,j)$.
\end{Lemma}
\begin{proof}
Let $b$ be the top half of $a$ and set $c=bb^{*}$. Then $c$ has  the same number of through strands as $a$, so $ (c,j)\in \tltrun$ for all $j\in\Z_m$. In particular: \[(a,i)=(c,i-j-\Phi(c,a))(a,j).\qedhere\]
\end{proof}

Recall from \autoref{L:UniqueFactor} that the top half of a diagram has no caps and the bottom half has no cups.  Let $B_l^n\subseteq C_l^n$ denote the set of diagrams without caps. 

\begin{Definition}
A pair $(a,b)\in B_l^n\times B_l^n$ is in \emph{vertical position} if $(a^*,0)(b,0)=(\Id_l,0)$ in $\mathbf{K}$. 
\end{Definition}

That is, the diagram $a^*b$ has no internal components and is equivalent to $\Id_l$.

\begin{Example}
The following pair of diagrams in $B^6_2$ are in vertical position:
\begin{gather*}
a=\;
\begin{tikzpicture}[anchorbase]
\draw[usual] (0.5,1) to[out=270,in=180] (0.75,0.75) to[out=0,in=270] (1,1);
\draw[usual] (-1,1) to[out=270,in=180] (-0.75,0.75) to[out=0,in=270] (-0.5,1);
\draw[usual] (0,0.5) to (0,1);
\draw[usual] (1.5,0.5) to (1.5,1);
\end{tikzpicture}
\;,\quad
b=\;
\begin{tikzpicture}[anchorbase,xscale=-1]
\draw[usual] (0.5,1) to[out=270,in=180] (0.75,0.75) to[out=0,in=270] (1,1);
\draw[usual] (-1,1) to[out=270,in=180] (-0.75,0.75) to[out=0,in=270] (-0.5,1);
\draw[usual] (0,0.5) to (0,1);
\draw[usual] (1.5,0.5) to (1.5,1);
\end{tikzpicture}
\;,\quad
a^{\ast}b
=\;
\begin{tikzpicture}[anchorbase, yscale=-1]
\draw[usual] (0.5,1) to[out=270,in=180] (0.75,0.75) to[out=0,in=270] (1,1);
\draw[usual] (-1,1) to[out=270,in=180] (-0.75,0.75) to[out=0,in=270] (-0.5,1);
\draw[usual] (0,0.5) to (0,1);
\draw[usual] (1.5,0.5) to (1.5,1);
\draw[usual] (0,1) to[out=90,in=0] (-0.25,1.25) to[out=180,in=90] (-0.5,1);
\draw[usual] (1.5,1) to[out=90,in=0] (1.25,1.25) to[out=180,in=90] (1,1);
\draw[usual] (0.5,1) to (0.5,1.5);
\draw[usual] (-1,1) to (-1,1.5);
\end{tikzpicture}
\;.\qedhere
\end{gather*}
\end{Example}

Let $\Gamma_l^n$ be the graph with vertex set $B_{l}^{n}$ and an edge between each pair of diagrams in vertical position.

\begin{Definition}
An \emph{outer cup} in a crossingless matching is a cup that is not separated from the bottom by another cup. A pair $(a,b)\in B_l^n\times B_l^n$ is a \emph{flip pair} if $ a$ is obtained from $ b$ by splitting an outer cup into a pair of through strands and connecting a pair of adjacent through strands into a cup.
\end{Definition}

\begin{Example}
The left diagram below has outer cups at $(2,3),(7,8)$ and adjacent through strands at $(5,8)$, so we can form a flip pair as follows:
\begin{gather*}
\begin{tikzpicture}[anchorbase]
\draw[usual] (0,0) to (0,0.5)node[above]{1};
\draw[usual] (0.5,0.5)node[above]{2} to[out=270,in=180] (0.75,0.3) to[out=0,in=270] (1,0.5)node[above]{3};
\draw[usual] (1.5,0) to (1.5,0.5)node[above]{4};
\draw[usual] (2,0) to (2,0.5)node[above]{5};
\draw[usual] (2.5,0.5)node[above]{6} to[out=270,in=180] (2.75,0.3) to[out=0,in=270] (3,0.5)node[above]{7};
\draw[usual] (3.5,0) to (3.5,0.5)node[above]{8};
\draw[ultra thick,spinach,dotted] (0.75,0) to (0.75,0.3);
\draw[ultra thick,spinach,dotted] (2,0.2) to (3.5,0.2);
\end{tikzpicture}
\ \longleftrightarrow \ 
\begin{tikzpicture}[anchorbase]
\draw[usual] (0,0) to (0,0.5)node[above]{1};
\draw[usual] (0.5,0) to (0.5,0.5)node[above]{2};
\draw[usual] (1,0) to (1,0.5)node[above]{3};
\draw[usual] (1.5,0) to (1.5,0.5)node[above]{4};
\draw[usual] (2,0.5)node[above]{5} to[out=270,in=180] (2.75,0.1) to[out=0,in=270] (3.5,0.5)node[above]{6};
\draw[usual] (2.5,0.5)node[above]{7} to[out=270,in=180] (2.75,0.3) to[out=0,in=270] (3,0.5)node[above]{8};
\end{tikzpicture}
\;.
\end{gather*}

The vertical dotted line indicates a splitting of an outer cup, and the horizontal dotted line indicates a closing of adjacent through strands. The resulting diagram now has an outer cup at $(5,6)$, and $(7,8)$ ceases to be an outer cup.
\end{Example}

Let $\Delta_l^n$ be the graph with vertex set $B_{l}^{n}$ and an edge between each flip pair.

\begin{Lemma}[{\cite[Lemma 4D.5 \& 4D.11]{Khovanov_2024}}]
We have that:
\begin{enumerate}
\item The graph $\Gamma_l^n$ is connected.
\item The graph $\Delta_l^n$ is connected for $l\ge 3$.
\end{enumerate}
\end{Lemma}

We now assume $l\ge 3$.

\begin{Lemma}\label{L:FlipVert}
\begin{enumerate}
\item[]

\item For $a,b,c\in B _{l}^{n}$, we have $ (ac^*,i)\approx_{l}(bc^{*},j)$ in $\tltrun$ for all $i,j\in\Z_m$.

\item For $k\le k'\le l$ and $a,b\in B_{k
'}^{n}$, we have $(aa^*,i)\approx_l(bb^*,j)$ in $\tltrun$  for all $i,j\in\Z_m$.

\end{enumerate}

\end{Lemma}
\begin{proof}
\begin{enumerate}
    
    \item First suppose that $(a,b)$ is in vertical position, so $(a^*,0)(b,0)=(\Id_l,0)$ in $\mathbf{K}$.  Then:
    \begin{align*}
    (a c^*,i)&=(aa^*b c^*,i)\\&=(aa^*,i-j-\Phi(aa^*, bc^*))(b c^*,j)\\&  \approx_l (bc^*, j).
    \end{align*}
    The claim then follows for arbitrary $(a,b)$ by the connectedness of $\Gamma^n_l$.

    \item First suppose that $(a,b)$ is a flip pair, and let $d\in B_{k'-2}^{n}$ be the diagram obtained from $a$ by closing up the pair of adjacent through strands involved in the flip pair. Then:
    \begin{align*}
        (aa^*,i)&\approx_l(dd^*,0)(aa^*,i)\\&=(dd^*,i+\Phi(dd^*,aa^*))\\
        &\approx_l (dd^*,j+\Phi(dd^*,bb^*))\\
        &=(dd^*,0)(bb^*,j)\\
        &\approx_l(bb^*,j).
    \end{align*}
    where the middle $\approx_l$ follows from \autoref{L:DiffBubLeftConn}. The claim then follows for arbitrary $(a,b)$ by the connectedness of $\Delta_{k'}^n$. \qedhere
\end{enumerate}
\end{proof}

Now recall the hook formulation for $m\tlmon[n]$.

\begin{Lemma}\label{L:TLLeftConn}
The monoid $\tltrun$ is left-connected.
\end{Lemma}
\begin{proof}
We first show that any diagram is left-connected to a diagram with $l$ through strands. Suppose $(a,i)\in\tltrun$ has a presentation $(a,i)=o^{i}u_{\alpha_r}u_{\alpha_{r-1}}\cdots u_{\alpha_1}$, with ${\alpha_s}\ne {\alpha_{s+1}}$ since we can push all internal components into $i$. Multiplying by a hook either preserves the number of through strands or decreases it by two, so there exists $1\le p\le r$ such that $(v,0)=u_{\alpha_p}\cdots u_{\alpha_1}$ has precisely $l$ through strands. Let $(v',0)=u_{\alpha_r}\cdots u_{\alpha_{p+1}}$, so $(a,0)=(v',0)(v,0)$. Since $(v,0)(v^*,0)(v,0)=(v,2c(l))$, we have:
\[(a,i)=(v',i-2c(l))(v,2c(l))=(v',i-2c(l))(v,0)(v^*,0)(v,0)\]
Since multiplication either preserves or reduces the number of through strands, the diagram $(v',i-2c(l))(v,0)(v^*,0)$ has between $k$ and $l$ through strands as needed, and so $(a,i)\approx_l(v,0)$.

It now suffices to show that every pair of diagrams $(v,j),(w,j')$ with $l$ through strands are left-connected. By \autoref{L:UniqueFactor}, we can express the diagrams as $(v,j)=(v_1v_2^*,j)$ and $(w,j')=(w_1w_2^*,j')$ with $v_1,v_2,w_1,w_2\in B_l^k$. Then by \autoref{L:FlipVert}:
\[(v,j)=(v_1v_2^*,j)\approx_l (v_2v_2^*,j)\approx_l (w_2w_2^*,j')\approx_l (w_1w_2^*,j')=(w,j').\qedhere\]
\end{proof}

\begin{Corollary}\label{L:TLRightConn}
The monoid $m\mathcal{TL}_{n}^{\le k}$ is right-connected.
\end{Corollary}
\begin{proof}
This follows from \autoref{L:TLLeftConn} by applying the diagrammatic anti-involution $*$.
\end{proof}

\begin{Lemma}\label{L:TLTrivialHom}
Suppose $\cchar\nmid m$. Then $\Hom(\tltrun,\K)=\{0\}$.
\end{Lemma}
\begin{proof}
Since every idempotent $ H$-cell in $\tltrun$ is either $\mathcal{C}_{m}$ or trivial, \autoref{T:Period} says that the period of every element divides $ m$. Therefore there exists $ N\in\N$ such that $(a,i)^{N}=(a,i)^{N+m}$ for all $(a,i)\in \tltrun$. Then for $f\in\Hom(\tltrun,\K)$, we have: \[Nf((a,i))=f((a,i)^N)=f((a,i)^{N+m})=(N+m)f((a,i)).\] It follows that $mf((a,i))=0$, and since $m\ne 0$ in $\K$,  $f$ is the zero map.
\end{proof}

Finally, we obtain that in good characteristic, the representation gap of $\tltrun$ is given by the minimal simple dimension.

\begin{Proposition}\label{P:TLGapIsMinDim}
Let $M$ be an $\tltrun$-representation with $\cchar\nmid m$. Then: \[\operatorname{gap}_{\K}(\tltrun)=\min\{ \dim_{\K}(L)\mid L\not\cong \onebt\text{ a simple }\tltrun\text{-representation} \}.\]
\end{Proposition}
\begin{proof}
Combining \autoref{L:TLNullConn}, \autoref{L:TLLeftConn}, and \autoref{L:TLRightConn}, we have that $\tltrun$ is well-connected. The claim then follows from \autoref{L:TLTrivialHom} and \autoref{T:WellConnect}.
\end{proof}

\subsection{Dimensions of simple representations}\label{S:TLSimples}

To find the minimal dimension of the non-trivial representations and hence the representation gap of $\tltrun$,  we need some results about the Temperley--Lieb algebra. 

\begin{Proposition}[{\cite[Theorem 6.7]{Graham1996}}]\label{P:TLAlgCellular}
The Temperley--Lieb algebra $\mathcal{TL}_{n}^{lin}(\delta)$ is cellular with cell datum $ (\mathcal P,\mathcal M,C,*) $ where:
\begin{enumerate}
\item $\mathcal P=\{ k\mid 0\le k\le n, \, k\equiv_{2}n \}$ with the reverse natural order.
\item $\mathcal M(k)=B_k^n$.
\item $c^{k}_{x,y}=xy^*$.
\end{enumerate}
The cell representations $W(k)$ are given by the $\K$-span of $B^n_k$, with the actions given by diagrammatic concatenation. The canonical bilinear forms $\langle-,-\rangle_k$ are given on $x,y\in B_k^n$ by:\[\langle x,y\rangle_k=
\begin{cases}
\delta^{\Phi(x^*,y)} &\text{if }x^*y=\Id_k,\\
0 &\text{else.}
\end{cases}
\]
\end{Proposition}

Similarly, the monoid algebra of $m\tlmon[n]$ has the structure of a sandwich cellular algebra, given by the unique factorisation from \autoref{L:UniqueFactor}.

\begin{Proposition}\label{P:CyclicTLSandwichDatum}
$(\K[m\tlmon[n]],m\tlmon[n])$ is an involutive sandwich pair with involutive sandwich cell datum:
\begin{enumerate}
\item $\mathcal{P}=\{ k\mid 0\le k\le n, \, k\equiv_{2}n \}$ with  the reverse natural order.
\item $\mathcal{T}(k)=B_k^n$.
\item $\mathscr{H}_k\cong\K\mathcal{C}_m$ with basis $B_\lambda=\{(\Id_k,0),\ldots,(\Id_k,m-1)\}$.
\item $c^{k}_{x,(\Id_k,i),y}=(xy^*,i)$.
\item $*$ given by the diagrammatic anti-involution.
\end{enumerate}
The cell representations $\Delta(k),\nabla(k)$ are given by the $\K$-span of $B^n_k\times \mathcal{C}_m$, with the actions given by diagrammatic concatenation. The bilinear maps $\phi^k$ are given by:\[\phi^k((x,i),(y,j))=
\begin{cases}
(\Id_k,i+j+\Phi(x^*,y)) &\text{if }x^*y=\Id_k,\\
0 &\text{else.}
\end{cases}\]
\end{Proposition}
\begin{proof}
By \autoref{P:BubCellStructure}, $m\tlmon[n]$ is regular, so $(\K[m\tlmon[n]],m\tlmon[n])$ is an involutive sandwich pair with respect to $*$ by \autoref{P:MonSandPair}. The same proposition provides the construction of the sandwich cell datum as follows. By \autoref{P:BubCellStructure}, the left (resp. right) cells have the same bottom (resp. top) half, so the top and bottom sets are the other halves $B^n_k$. The sandwich cellular basis elements $c^{k}_{x,(\Id_k,i),y}$ are given by diagrammatic concatenation, so $c^{k}_{x,(\Id_k,i),y}=(x,0)(\Id_k,i)(y^*,0)=(xy^*,i)$. The cell representations follow by construction.
\end{proof}

Suppose $\K$ contains a primitive $m$-th root of unity, and let $\zeta_m$ be an $m$-th root of unity with corresponding simple $\mathcal{C}_m$-representation $K$. Recall the notation that $L(k,K)$ denotes the simple $m\tlmon[n]$-representation  associated with the $k$-th cell and $K$. Let $L(k)$ denote the $k$-th simple representation for $\mathcal{TL}_{n}^{lin}(\zeta_m)$ as indexed by $\mathcal{P}$.

\begin{Theorem}\label{T:SimpleDims}
Suppose $\K$ contains a primitive $m$-th root of unity. Then: \[\dim_\K L(k,K)=\dim_\K L(k).\]
\end{Theorem}
\begin{proof}
By cellular algebra theory, we have $L(k)=W(k)/R(k)$, where $W(k)$ is the $k$-th cell representation of $\tlalg[n]{\zeta_m}$ with basis $B_k^n$, and $R(k)$ is the kernel of the canonical bilinear form $\langle -,-\rangle_k$ as in \autoref{P:TLAlgCellular}. By \autoref{T:SandwichHReduce}, we have $L(k,K)=\Delta(k,K)/R(k,K)$, where $\Delta(k,K)=\Delta(k)\otimes_{\K\mathcal{C}_{m}}K$, and $R(k,K)$ is the kernel of the map $\Delta(k,K) \to\Hom_\K(\nabla(k),\mathcal{C}_m)\otimes_{\K\mathcal{C}_{m}}K$. 

Now $\Delta(k,K)$ is the quotient of $\Delta(k)\otimes_{\K} K$ by the $\K$-span of elements of the form $((x,i)\cdot g)\otimes 1 - (x,i)\otimes (g\cdot 1)$, where $g\in\mathcal{C}_m$ acts on $\Delta(k)$ by adding an internal component and on $K$ by multiplication by $\zeta_m$. It follows that basis elements of $\Delta(k,K)$ satisfy: \[(x,i)\otimes 1=((x,0)\cdot g^i)\otimes 1=(x,0)\otimes(g^i\cdot 1)=(x,0)\otimes \zeta_m^i=\zeta_m^i((x,0)\otimes 1).\]
Since $\Delta(k)\otimes_{\K} K\cong \Delta(k)$ as vector spaces, it follows that  $\Delta(k)$ under the identification $(x,i)=\zeta_m^i(x,0)$ is isomorphic to $\Delta(k,K)$. In particular, $B_k^n$ is a basis for the former, so $\dim_\K(\Delta(k,K))=|B_k^n|=\dim_\K(W(k))$.

Similarly identifying $\Hom_\K(\nabla(k),\mathcal{C}_m)\otimes_{\K\mathcal{C}_{m}}K$ as a quotient of $\Hom_\K(\nabla(k),\mathcal{C}_m)$, we have $\phi^k( -, (x,i) )= \zeta_m^i\phi^k( -, (x,0) )$. Let $(y,i)$ be a basis element of $\nabla(k)$. Then:
\[
\phi^k((y,i),(x,0))=\phi^k((y,0),(x,i))=\zeta_m^i\phi^k( (y,0), (x,0) ).
\]
Let $\sum_\alpha \mu_\alpha x_\alpha\in W(k)$ with $\mu_\alpha\in\K$ and $x_\alpha\in B_k^n$. Then:
\begin{align*}
\phi^k\left((y,i),\sum_\alpha \mu_\alpha (x_\alpha,0)\right)&=\zeta_m^i\sum_\alpha \mu_\alpha \phi^k\left((y,0), (x_\alpha,0)\right)\\
&=\zeta_m^i\sum_\alpha \mu_\alpha \langle y,x_\alpha\rangle_k \phi^k\left((y,-\Phi(y,x_\alpha)), (x_\alpha,0)\right) \\
&=\zeta_m^i\sum_\alpha \mu_\alpha \langle y,x_\alpha\rangle_k \\
&=\zeta_m^i\left\langle y,\sum_\alpha \mu_\alpha x_\alpha \right\rangle_k.
\end{align*}
It follows that $\sum_\alpha \mu_\alpha x_\alpha\in R(k)$ if and only if $\sum_\alpha \mu_\alpha (x_\alpha,0)\in R(k,K)$, so $\dim_\K(R(k))=\dim_\K(R(k,K))$. The stated equality now follows.\end{proof}

\begin{Example}
Fields containing a primitive $m$-th root of unity include algebraically closed fields with characteristic not dividing $m$ and finite fields $\F_{p^k}$ with $m\mid p^k-1$.
\end{Example}

It remains to compute the dimensions, for which we need the following definitions.

\begin{Definition}
The \emph{Chebyshev polynomials} (of the second kind and normalised) in $\K[X]$ are given by the recurrence relations:
\begin{align*}
U_0(X)&=1,\\
U_1(X)&=X,\\
U_k(X)&=XU_{k-1}(X)-U_{k-2}(X).
\end{align*}
The \emph{quantum characteristic} of $\delta\in\K$ is the minimal $l\in\N$ such that $U_{l+1}(\delta)=0$, with $l=\infty$ if $U_k(\delta)$ never vanishes. Let $p=\cchar[\K]$ (with $p=\infty$ if $\cchar[\K]=0$), and $\nu_p$ denote $p$-adic valuation. For $x\in\N$, we define $\nu_{l,p}(x)=\nu_{p}\left( \frac{x}{l} \right)$ {if }$x\equiv_l 0$, and $\nu_{l,p}(x)=0$ otherwise.
The \emph{$l,p$-adic expansion} of $x$ is given by:
\[x=\sum_{i=1}^{\infty}lp^{i-1}x_{i}+x_{0}=:[\dots,x_{1},x_{0}]\]
for $x_{i>0}\in\{ 0,\dots,p-1 \}$ and $x_{0}\in \{ 0,\ldots,l-1 \}$.
We also write:
\begin{enumerate}
\item $x\lhd y$ if $[\dots,x_{1},x_{0}]$ is digit-wise less than or equal to $[\dots,y_{1},y_{0}]$.
\item $x\lhd ' y$ if $x\lhd y$, $\nu_{l,p}(x)=\nu_{l,p}(y)$, and the $\nu_{l,p}(x)$-th digit of $x,y$ agree.
\end{enumerate}
And lastly:
\[e_{n,k}:=\begin{cases}
1 & \text{if }n\equiv_2 k,\,\nu_{l,p}(k)=\nu_{l,p}\left( \frac{n+k}{2} \right),\,k\lhd'\frac{n+k}{2}. \\
-1 & \text{if }n\equiv_2 k,\,\nu_{l,p}(k)<\nu_{l,p}\left( \frac{n+k}{2} \right),\,k\lhd\frac{n+k}{2}-1. \\
0 & \text{else.}
\end{cases}\]
\end{Definition}

\begin{Theorem}[{\cite[Corollary 9.3]{Spencer_2023}}]\label{P:TLAlgDim}
Let $L(k)=W(k)/R(k)$ be the $k$-th simple $\mathcal{TL}_{n}^{lin}(\delta)$-representation. We have:\[\dim_{\K}(L(k))=\sum_{r=0}^{c(k)}e_{n-2r+1,k+1}\left( \frac{n-2r+1}{n-r+1}\binom{n}{r} \right).\]
In particular, $\dim_{\K}(L(k))=1$ for $k\in\{0,1,n\}$.
\end{Theorem}

We have the desired formula for the dimensions of the simple representations, however it is difficult to compute due to the dependence on the base field. Nevertheless, in certain cases we have a easily stated lower bound.

\begin{Corollary}\label{C:DimLowerBound}
Suppose $\K$ contains a primitive $m$-th root of unity. If $\cchar[\K]=0$ or $k\ge\sqrt{n+2}-2$, then:
\[\dim_{\K}(L(k,K))\ge\frac{1}{(n-c(k)+1)(n-c(k)+2)}\binom{n}{c(k)}.
\]
\end{Corollary}
\begin{proof}
This follows from \autoref{T:SimpleDims}, \autoref{P:TLAlgDim} and \cite[Proposition 9.4 \& 9.5]{Spencer_2023}.
\end{proof}

\begin{Remark}

We'd like to point out that there are alternative approaches to proving \autoref{T:SimpleDims} without using sandwich cellular theory. For example, one can identify each $\tlalg[n]{\zeta_m^i}$ as a quotient algebra of $\K[m\tlmon[n]]$, show that the induced representations are non-isomorphic, and deduce that they form a complete set of simple $m\tlmon[n]$-representations by $H$-reduction and the cellularity of $\tlalg[n]{\zeta_m^i}$. But in a general setting, the step of showing non-isomorphism is usually non-trivial, whereas the approach in \autoref{T:SimpleDims} can be generalised to the problem of computing determinants of \emph{sandwich matrices} associated with the bilinear maps  (see \cite[Section 2D]{Tubbenhauer_2024}). For future work, one can use this to obtain exact bounds when $\K$ does not contain a primitive $m$-th root of unity. Moreover, the sandwich cellular approach also generalises to diagram monoids with arbitrary maximal subgroups and non-involutive diagram monoids (see \cite[Section 4]{Tubbenhauer_2024} for the transformation monoid and symmetric monoids).
\end{Remark}

\subsection{Representation gap}\label{S:TLRepGap}

Here, we bring all the results above together to establish the representation gap of the cyclic Temperley--Lieb monoid, as well as related measures. We retain the notation from the previous subsection.

\begin{Theorem}\label{T:GapIsBoundedByDim}

Suppose $\cchar[\K]\nmid m$.  Then:
\begin{align*}
\operatorname{gap}_{\K}(m\mathcal{TL}_{n}^{k\le l})&\ge\min\{ \dim_\K(L(j,K))\mid j\in \{ k,k+2,\ldots,l \} \}
\end{align*}

\end{Theorem}
\begin{proof}
Since $\cchar[\K]\nmid m$, the algebraic closure of $\K$ contains a primitive $m$-th root of unity. The inequality then follows from \autoref{T:FieldExtGap}, \autoref{P:TLGapIsMinDim} and \autoref{T:SimpleDims}.
\end{proof}

Recall that the Big Theta notation $\Theta(f)$ means asymptotically bounded above and below by $f$.

\begin{Theorem}\label{T:TLRepGap}

Suppose $\cchar[\K]\nmid m$, and $\mathbb L$ is an arbitrary field. Let $\sqrt{n+2}-2\le k\le l\le 2\sqrt{n}$.  Then:
\begin{align*}\operatorname{gap}_{\K}(m\mathcal{TL}_{n}^{k\le l}) & \ge \frac{4}{(n+2\sqrt{ n }+2)(n+2\sqrt{ n }+4)}\binom{n}{\lfloor\frac{n}{2}-\sqrt{ n }\rfloor}\in \Theta(2^{n}n^{-5/2}), \\\operatorname{ssgap}_{\mathbb{L}}(m\mathcal{TL}_{n}^{k\le l}) & \ge \frac{4\sqrt{n}-2}{n+2\sqrt{n}+1}\binom{n}{\lfloor\frac{n}{2}-\sqrt{ n }\rfloor}\in\Theta(2^{n}n^{-1}).\end{align*}

\end{Theorem}
\begin{proof}
Let $k\le j\le l$. The bound for the representation gap comes from applying \autoref{T:GapIsBoundedByDim} to the formula
$\frac{1}{(n-c(j)+1)(n-c(j)+2)}\binom{n}{c(j)}$ in \autoref{C:DimLowerBound}. This function has its minimum at $j=\lfloor 2\sqrt{n}\rfloor$, plugging which into the formula yields the stated bound. 

The bound for the semisimple representation gap comes from the formula  $\frac{j+1}{n-c(j)+1}\binom{n}{c(j)}$ in \autoref{P:SSDimFormula}. This function has its maximum at $\lfloor\frac{\sqrt{8n+17}-1}{4}\rfloor$, and is monotonically increasing to its left and monotonically decreasing to its right. Since $\sqrt{n+2}-2> \frac{\sqrt{8n+17}-1}{4}$ for large enough $n$, the semisimple dimension is minimised at $j=\lfloor2\sqrt{n}\rfloor$, giving the stated bound.

Both asymptotic formulae are obtained from applying Stirling's approximation for binomial coefficients.
\end{proof}

We can then derive the associated ratios. Recall that the Big Omega notation $\Omega(f)$ means asymptotically bounded below by $f$.

\begin{Corollary}\label{C:TLRatios}
Suppose $\cchar[\K]\nmid m$, and $\mathbb L$ is an arbitrary field. Let $\sqrt{n+2}-2\le k\le l\le 2\sqrt{n}$.  Then:
\begin{align*}\operatorname{gapr}_{\K}(m\mathcal{TL}_{n}^{k\le l}) & \in \Omega(n^{-7/4}), \\\operatorname{ssgapr}_{\mathbb{L}}(m\mathcal{TL}_{n}^{k\le l}) & \in\Omega(n^{-1/4}).\end{align*}
\end{Corollary}
\begin{proof}
By  \autoref{P:BubCellStructure}, we have $|\mathcal{J}_j|=m\left(\frac{j+1}{n-c(j)+1}\binom{n}{c(j)}\right)^2$. By the proof of \autoref{T:TLRepGap}, this is maximised at $j=\lceil \sqrt{n+2}\rceil+2$, so:
\[|\tltrun|\le m\left(\frac{2 \sqrt{n+2}+8}{n+ \sqrt{n+2}+3}\binom{n}{\lceil \frac{n-\sqrt{n+2}}{2}\rceil-1}\right)^2\left(\frac{2\sqrt{n}-\sqrt{n+2}-2}{2}\right)\]

Applying asymptotics similarly to \autoref{T:TLRepGap} and plugging the formula into the corresponding ratios yields the stated bounds.
\end{proof}

As per \autoref{C:DimLowerBound}, in characteristic zero we can preserve the bound on the representation gap without needing to truncate $m\tlmon[n]$ from below. In this case, we also have a finer bound on the faithfulness.

\begin{Theorem}
Suppose $\cchar[\K]=0$, and $\mathbb L$ is an arbitrary field. Let $0\le l\le 2\sqrt{n}$.  Then:
\begin{align*}\operatorname{gap}_{\K}(m\mathcal{TL}_{n}^{\le l}) & \ge \frac{4}{(n+2\sqrt{ n }+2)(n+2\sqrt{ n }+4)}\binom{n}{\lfloor\frac{n}{2}-\sqrt{ n }\rfloor}\in \Theta(2^{n}n^{-5/2}), \\\operatorname{ssgap}_{\mathbb{L}}(m\mathcal{TL}_{n}^{\le l}) & \ge \frac{1}{n}\binom{n}{\lfloor\frac{n}{2}\rfloor}\in\Theta(2^{n}n^{-3/2}),\\
\operatorname{faith}_{\K}(m\mathcal{TL}_{n}^{\le l})  &\ge \frac{6}{n+4}\binom{n}{\lfloor\frac{n}{2}\rfloor-1}\in\Theta(2^{n}n^{-3/2}).
\end{align*}

\end{Theorem}
\begin{proof}
The bound for the representation gap is unchanged from \autoref{T:TLRepGap}.

For the semisimple representation gap, we now need to compare $j=0,1$ and $j=\lfloor 2n\rfloor$. Simple algebra shows the former case has smaller semisimple representation dimension, yielding the stated bound.

For the faithfulness, by \autoref{L:FaithEmbed}, we only need to bound $\faith[\K]{m\tlmon[n]^{\le 2}}$ since there is an embedding $m\tlmon[n]^{\le 2}\hookrightarrow m\tlmon[n]^{\le l}$ (adding a strand to $m\tlmon[n]^{\le 2}$ in the case when $n$ is odd). By \autoref{P:TLGapIsMinDim}, faithful $m\tlmon[n]^{\le 2}$ representations cannot be an extension of only $\onebt$'s, so it must contain $L(2,K)$ for some simple $\mathcal{C}_m$-representation $K$. By \cite[Example 4B.9]{Khovanov_2024}, we have $\dim_\K(L(2,K))=\operatorname{ssdim}_\K(L(2,K))$ in characteristic zero, so computing the latter using \autoref{P:SSDimFormula} yields the stated bound.
\end{proof}

\begin{Remark}
The associated ratios in the characteristic zero case can be computed analogously to \autoref{C:TLRatios}. In particular, the gap ratio remains unchanged.
\end{Remark}

\begin{Remark}
We chose $l\le 2\sqrt{n}$ above to keep the simple dimensions, and hence the representation gap, large. However, fewer through strands in the monoid destroys information quicker during multiplication, which makes the key space smaller for the attacker. To balance the level of security between the representation gap and the size of the key space, one can use a larger truncation by choosing an $l$ that is closer to $n$. The method for finding the bounds in this case is analogous to \autoref{T:TLRepGap}, and is computed for the ordinary Temperley--Lieb monoid in \cite{Khovanov_2024}.
\end{Remark}

\section{Other cyclic diagram monoids}\label{Chapter4}

We now discuss cyclic versions of other planar diagram monoids listed in the introduction. We can similarly define these monoids as endomorphism monoids of monoidal categories as in the Temperley--Lieb case, but for brevity, we define them as diagrams of set partitions as per \autoref{R:PartitionDiagram}, and focus on their representation gaps. We also retain our terminologies and notations from before, e.g. the ordered pair notation for diagrams.

\subsection{Cyclic planar partition monoids}

The \emph{$m$-cyclic planar partition monoid} $m\ppamon$ consists of planar diagrams of arbitrary partitions of $2n$ points $\{1,\ldots, n,1',\ldots,n'\}$, with internal components counted modulo $m$.

\begin{Example}
The following is a diagram in $3\ppamon[4]$ with two internal components, where strands connect points in the same set:
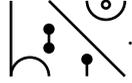
\[
\begin{tikzpicture}[anchorbase]
\draw[usual] (0,0) to[out=90,in=180] (0.25,0.25) to[out=0,in=90] (0.5,0);
\draw[usual] (1,1) to[out=270,in=180] (1.25,0.75) to[out=0,in=270] (1.5,1);
\draw[usual] (0,0) to (0,1);
\draw[usual] (1.5,0) to (0.5,1);
\draw[usual, dot] (.5,.5) to (.5,.35);
\draw[usual, dot] (.5,.5) to (.5,.65);
\draw[usual] (1.25,0.95) circle [radius=0.5mm];
\draw[usual,dot] (1,0) to (1,.25);
\end{tikzpicture}.
\]
The dot extending out from a bottom point indicates a singleton, which we call \emph{stops}.
\end{Example}

Note that the diagram of a partition is not unique with regards to the strands connecting points, but they are well-defined under multiplication.

\begin{Proposition}
The monoid $m\ppamon$ has presentation:
\begin{equation*}
    m\ppamon= \scaleleftright[1.75ex]{<} {c,p_1,p_{\frac{3}{2}},p_2,\ldots,p_{n-\frac{1}{2}}, p_n\, \vrule width 1pt \, \begin{matrix}cp_i=p_ic,\\ p_i^2=cp_i, \text{ for }i\in\Z,\\ p_i^2=p_i, \text{ for }i\in\Z+\frac{1}{2},\\ p_ip_{i\pm\frac{1}{2}}p_i=p_i,\\ p_ip_j=p_jp_i, \text{ for } |i-j|>\frac{1}{2},\\ c^m=1 \end{matrix}} {>}.
\end{equation*}
\end{Proposition}
\begin{proof}
Analogous to the Temperley--Lieb case in \autoref{T:TLPres} using e.g. \cite[Theorem 1.11]{HALVERSON2005869}.
\end{proof}

The generator $c$ is the identity diagram with one internal component, whereas $p_i$ and $p_{i+\frac{1}{2}}$ are given by:
\[
p_i=\begin{tikzpicture}[anchorbase, scale=1.5, every node/.style={scale=0.75}]
\draw[usual] (0,0) to (0,1)node[above]{$i-1$};
\draw[usual] (1.5,0) to (1.5,1)node[above]{$i+2$};
\node at (-0.5,0.5) {$\cdots$};
\node at (2,0.5) {$\cdots$};
\draw[usual] (-1,0) to (-1,1)node[above]{$1$};
\draw[usual] (2.5,0) to (2.5,1)node[above]{$n$};
\draw[usual, dot] (.5,0) to (.5,.25);
\draw[usual, dot] (.5,1)node[above]{$i$} to (.5,.75);
\draw[usual] (1,0) to (1,1)node[above]{$i+1$};
\end{tikzpicture}\quad,\quad
p_{i+\frac{1}{2}}=\begin{tikzpicture}[anchorbase, scale=1.5, every node/.style={scale=0.75}]
\draw[usual] (0,0) to (0,1)node[above]{$i-1$};
\draw[usual] (1.5,0) to (1.5,1)node[above]{$i+2$};
\node at (-0.5,0.5) {$\cdots$};
\node at (2,0.5) {$\cdots$};
\draw[usual] (-1,0) to (-1,1)node[above]{$1$};
\draw[usual] (2.5,0) to (2.5,1)node[above]{$n$};
\draw[usual] (1,0) to (1,1)node[above]{$i+1$};
\draw[usual] (.5,0) to (.5,1)node[above]{$i$};
\draw[usual] (0.5,1) to[out=270,in=180] (0.75,0.75) to[out=0,in=270] (1,1);
\draw[usual] (0.5,0) to[out=90,in=180] (0.75,0.25) to[out=0,in=90] (1,0);
\end{tikzpicture}.
\]

By e.g. \cite[(1.5)]{HALVERSON2005869}, we have a monoid isomorphism $\ppamon\cong\tlmon[2n]$, but more generally we have the following.

\begin{Proposition}\label{P:PaEmbedsTL}
There is an embedding of monoids:
\[m\ppamon[n]\hookrightarrow 2m\tlmon[2n].\]
\end{Proposition}
\begin{proof}
We define a map:
\begin{align*}
p_i &\mapsto ou_{2i-1},\\
p_{i+\frac{1}{2}} &\mapsto o^{2m-1}u_{2i},\\
c &\mapsto o^2.
\end{align*}
One can verify that the relations are satisfied in $2m\tlmon[2n]$, so we get a homomorphism. For injectivity, suppose $(a,i)\in 2m\tlmon[2n]$ has a presentation $(a,i)=o^{i}u_{\alpha_r}\cdots u_{\alpha_1}$, then $p_{\frac{\alpha_r+1}{2}}\cdots p_\frac{\alpha_1+1}{2}$ also has skeleton $a$ under the homomorphism. Since the image contains all $o^2,\ldots,o^{2m-2}$, it follows that the image also contains either all diagrams with skeleton $a$ and odd numbers of internal components or  all diagrams with skeleton $a$ and even numbers of internal components. Since $|m\ppamon|=|2m\tlmon[2n]|/2$, the homomorphism is injective.
\end{proof}

Note that the above also shows that each $J$-cell of $2m\tlmon[2n]$ restricts to a $J$-cell of $m\ppamon$. Indeed, an analogous statement to \autoref{C:BubCell} holds for $m\ppamon$, which implies it has the same cell structure as $\ppamon\cong\tlmon[2n]$ but with each $H$-cell inflated. In particular, its $J$-cells are ordered by the number of through strands. We then have the following.

\begin{Corollary}\label{C:pPaLessTL}
Let $\mathcal J$ be a $J$-cell of $2m\tlmon[2n]$, and $L_\mathcal{J}^{2m\tlmon[2n]}$ (resp. $L_\mathcal{J}^{m\ppamon}$)  be an associated simple $2m\tlmon[2n]$-representation (resp. $m\ppamon$-representation). Then:
\[\dim_\K(L_\mathcal{J}^{m\ppamon})\le\dim_\K(L_\mathcal{J}^{2m\tlmon[2n]}). \]
\end{Corollary}
\begin{proof}
This follows from the discussion above and \cite[Theorem 3D.2]{Khovanov_2024}.
\end{proof}

Let $m\ppamon^{k\le l}$ denote the truncation of $m\ppamon$ to between $k$ and $l$ through strands. 

\begin{Corollary}
We have:
\[\faith[\K]{m\ppamon^{k\le l}}\le \faith[\K]{2m\tlmon[2n]^{2k\le 2l}}.\]
\end{Corollary}
\begin{proof}
This follows from the embedding and \autoref{L:FaithEmbed}.
\end{proof}

Moreover, we see that the representation gap of cyclic planar partition monoids is also no better than cyclic Temperley--Lieb monoids. 
\begin{Theorem}
Suppose $\cchar\nmid 2m$. Then:
\[\gap[\K]{m\ppamon^{k\le l}}\le \gap[\K]{2m\tlmon[2n]^{2k\le 2l}}.\]
\end{Theorem}
\begin{proof}
By \autoref{L:GapBound}, \autoref{C:pPaLessTL}. and \autoref{P:TLGapIsMinDim}:
\begin{align*}
\gap[\K]{m\ppamon^{k\le l}}&\le\min \{ \dim_{\K}(L_{j})\mid L_{j}\not\cong \onebt\text{ a simple }{m\ppamon^{k\le l}}\text{-representation} \}\\
&\le\min \{ \dim_{\K}(L_{j})\mid L_{j}\not\cong \onebt\text{ a simple }{2m\tlmon[2n]^{2k\le 2l}}\text{-representation} \}\\
&=\gap[\K]{2m\tlmon[2n]^{2k\le 2l}}.\qedhere
\end{align*}
\end{proof}

\subsection{Cyclic Motzkin monoids}

The \emph{$m$-cyclic Motzkin monoid} $m\momon$ consists of planar diagrams of partitions of $2n$ points $\{1,\ldots, n,1',\ldots,n'\}$ into sets of either size one or size two, with internal components counted modulo $m$. 

\begin{Example}
The following is a diagram in $3\momon[4]$ with two internal components:
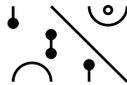
\[
\begin{tikzpicture}[anchorbase]
\draw[usual] (0,0) to[out=90,in=180] (0.25,0.25) to[out=0,in=90] (0.5,0);
\draw[usual] (1,1) to[out=270,in=180] (1.25,0.75) to[out=0,in=270] (1.5,1);
\draw[usual] (1.5,0) to (0.5,1);
\draw[usual, dot] (.5,.5) to (.5,.35);
\draw[usual, dot] (.5,.5) to (.5,.65);
\draw[usual] (1.25,0.95) circle [radius=0.5mm];
\draw[usual,dot] (1,0) to (1,.25);
\draw[usual,dot] (0,1) to (0,.75);
\end{tikzpicture}.\qedhere
\]
\end{Example}

The Motzkin algebra $\momon^{lin}(\delta)$ (see \cite{BENKART2014473}) is constructed similarly as the Temperley--Lieb algebra is to the Temperley--Lieb monoids,  but we do not know its simple dimensions, so we have the following but not the representation gap. Let $m\momon^{\le k}$ denote the truncation of $m\momon$ to at most $k$ through strands.

\begin{Theorem}
\begin{enumerate}
    \item[]
    \item The monoid $m\momon$ has the same cell structure as $\momon$ with each $H$-cell inflated to size $m$.
    \item Each simple representations $L$ of $m\momon$ is indexed by $(j,K)$ where $j\in\{n,n-2,\ldots\}$ and $K$ is a simple $\mathcal{C}_m$-representation.
    \item If $\cchar\nmid m$, then $m\momon^{\le k}$ is well-connected, and:\[\operatorname{gap}_{\K}(m\momon^{\le k})=\min\{ \dim_{\K}(L)\mid L\not\cong \onebt\text{ a simple }m\momon^{\le k}\text{-representation} \}.\]
    \item The monoid algebra $\K[m\momon]$ is involutive sandwich cellular with involutive sandwich cell datum given by \autoref{P:MonSandPair}.
    \item If $\K$ contains a primitive root of unity, $\zeta_m$ is an $m$-th root of unity with corresponding simple $\mathcal{C}_m$-representation $K$, and $L(j)$ is the $j$-th simple representation for $\momon^{lin}(\zeta_m)$, then:\[\dim_\K L(j,K)=\dim_\K L(j).\]
\end{enumerate}
\end{Theorem}
\begin{proof}
Analogous to the Temperley--Lieb case using \cite{2025arXiv251006707A}.
\end{proof}

\subsection{Cyclic planar rook monoids}

The \emph{$m$-cyclic planar rook monoid} $m\promon$ consists of planar diagrams of partitions of $2n$ points $\{1,\ldots, n,1',\ldots,n'\}$ into sets of either size one or size two, but with the latter consisting of exactly one element from each of $\{1,\ldots, n\}$ and $\{1',\ldots,n'\}$, with internal components counted modulo $m$.

\begin{Example}
The following is a diagram in $3\promon[4]$ with one internal component:
\[
\begin{tikzpicture}[anchorbase]
\draw[usual] (0,0) to (0,1);
\draw[usual] (1.5,0) to (0.5,1);
\draw[usual, dot] (.25,.5) to (.25,.35);
\draw[usual, dot] (.25,.5) to (.25,.65);
\draw[usual,dot] (.5,0) to (.5,.25);
\draw[usual,dot] (1,0) to (1,.25);
\draw[usual,dot] (1,1) to (1,.75);
\draw[usual,dot] (1.5,1) to (1.5,.75);
\end{tikzpicture}.\qedhere
\]
\end{Example}

\begin{Proposition}
The monoid algebra $\K[m\promon]$ is semisimple if and only if $\cchar\nmid m$.
\end{Proposition}
\begin{proof}
By \cite[Corollary 9.4]{Steinberg2016}, it suffices to show that $m\promon$ is an inverse monoid (i.e. for every element $a$, there exists a unique element $a^*$ such that $aa^*a=a$ and $a^*aa^*=a$), which is equivalent to the monoid being regular and having commuting idempotents by \cite[Theorem 3.2]{Steinberg2016}. 

By an analogous statement to \autoref{C:BubCell}, $m\promon$ has the same cell structure as $\promon$ but with each $H$-cell inflated to size $m$, in particular having the same idempotent $H$-cells, so it is regular by e.g. \cite[Proposition 4F.4]{Khovanov_2024}. 

We claim that the idempotents of $m\promon$ contain no diagonal through strands, that is, their corresponding partitions contain no sets of the form $\{i,j'\}$ with $i\ne j$. Indeed, we would otherwise need $\{i,i'\}$ and $\{j,j'\}$ to be in the partition, a clear contradiction. So for each $i=1,\ldots,n$, an idempotent either contains  $\{i,i'\}$ or both $\{i\},\{i'\}$, that is, either a straight through strand or  two stops. Since:
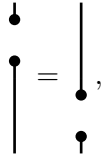
\[
\begin{tikzpicture}[anchorbase]
\draw[usual] (0,0) to (0,1);
\draw[usual,dot] (0,1) to (0,1.25);
\draw[usual,dot] (0,2) to (0,1.75);
\end{tikzpicture}\;=\;
\begin{tikzpicture}[anchorbase]
\draw[usual] (0,1) to (0,2);
\draw[usual,dot] (0,0) to (0,.25);
\draw[usual,dot] (0,1) to (0,.75);
\end{tikzpicture}\;
,
\]
it follows that idempotents of $m\promon$ commute.
\end{proof}

 By \autoref{P:ssdim}, the semisimplicity in good characteristic means the representation gap is easy to obtain. Let $m\promon^{k\le l}$ denote the  truncation  of $m\promon$ to between $k$ and $l$ through strands.

\begin{Theorem}
Suppose $\cchar\nmid m$. Then:
\[\gap[\K]{m\promon^{k\le l}}\ge\min\left\{\binom{n}{k},\binom{n}{l}\right\}.\]
In particular:
\[\gap[\K]{m\promon^{\lfloor\frac{n}{2}-\sqrt{ n }\rfloor\le \lceil\frac{n}{2}+\sqrt{ n }\rceil}}\ge\binom{n}{\lfloor\frac{n}{2}-\sqrt{ n }\rfloor}\in\Theta(2^nn^{-1/2}).\]
\end{Theorem}
\begin{proof}
Similar to \cite[Theorem 4F.15]{Khovanov_2024} and omitted.
\end{proof}


\bibliographystyle{alphaurl}
\bibliography{references}

\end{document}